\documentclass[a4paper,fleqn]{cas-dc}

\usepackage{float}
\usepackage[numbers]{natbib}

\usepackage{booktabs}
\usepackage[singlelinecheck=false]{caption}
\usepackage{algorithm}
\usepackage{graphicx}
\usepackage{lipsum}
\usepackage{algpseudocode}
\usepackage{amsmath}
\usepackage{amssymb}
\usepackage{booktabs}
\usepackage{arydshln}
\usepackage{lipsum}
\usepackage{natbib}
\usepackage{comment}
\usepackage{color}
\usepackage{cases}
\usepackage{adjustbox}
\usepackage{seqsplit}
\usepackage{mwe}
\usepackage{multirow}
\usepackage{epstopdf}
\usepackage{csquotes}
\usepackage{stfloats}
\usepackage{amsthm}
\usepackage{adjustbox}
\usepackage{colortbl}
\usepackage{lineno,hyperref}
\usepackage{babel}
\usepackage{multicol}
\usepackage{multirow}
\usepackage{siunitx}
\usepackage{float} 
\usepackage{enumitem}

\newcommand*{\rom}[1]{\expandafter\@slowromancap\romannumeral #1@}
\newtheorem{theorem}{Theorem}[section]

\newtheorem{proposition}[theorem]{Proposition}

\def\tsc#1{\csdef{#1}{\textsc{\lowercase{#1}}\xspace}}
\tsc{WGM}
\tsc{QE}
\tsc{EP}
\tsc{PMS}
\tsc{BEC}
\tsc{DE}
\begin{document}  \sloppy
% \let\WriteBookmarks\relax
% \def\floatpagepagefraction{1}
% \def\textpagefraction{.001}
% \shorttitle{}
% \shortauthors{Meysam Mahjoob}
% \author{Meysam Mahjoob}
% \author{Seyed Sajjad Fazeli}
% \author{Soodabeh Milanloui}
% \author{Leyla Saadat Tavasoli}
% \title [mode = title]{A Modified Adaptive Genetic Algorithm for Multi-product Multi-period Inventory Routing Problem}       
% \author{Meysam Mahjoob}

% \begin{abstract}
% Recent developments in urbanization and e-commerce have pushed businesses to deploy efficient systems to decrease their supply chain cost. Vendor Managed Inventory (VMI) is one of the most widely used strategies to effectively manage supply chains with multiple parties. VMI implementation asks for solving the Inventory Routing Problem (IRP). This study considers a multi-product multi-period inventory routing problem, including a supplier, set of customers, and a fleet of heterogeneous vehicles. Due to the complex nature of the IRP, we developed a Modified Adaptive Genetic Algorithm (MAGA) to solve a variety of instances efficiently. As a benchmark, we considered the results obtained by Cplex software and an efficient heuristic from the literature. Our approach showed significant improvement comparing the other two methods.
% \end{abstract}
% \begin{keywords}
% Inventory Routing Problem \sep Genetic Algorithm \sep Vendor Managed Inventory \sep Supply Chain Management \sep Adaptive Heuristic
% \end{keywords}

\let\WriteBookmarks\relax
\def\floatpagepagefraction{1}
\def\textpagefraction{.001}
\shorttitle{}
% \shortauthors{Meysam Mahjoob et~al.}

\title [mode = title]{Green Supply Chain Network Design with Emphasis on Inventory Decisions}                      
% \tnotemark[1,2]

% \tnotetext[1]{This document is the results of the research
%   project funded by the National Science Foundation.}

% \tnotetext[2]{The second title footnote which is a longer text matter
%   to fill through the whole text width and overflow into
%   another line in the footnotes area of the first page.}

\author[1]{Meysam Mahjoob}[
% type=editor,
%                         auid=000,bioid=1,
%                         prefix=Sir,
%                         role=Researcher,
%                         orcid=0000-0001-7511-2910
                        ]

% \fnmark[1]
\ead{mahjoob_m@alumni.ut.ac.ir}
% \ead[url]{www.cvr.cc, cvr@sayahna.org}

% \credit{Conceptualization of this study, Methodology, Software}

\address[1]{Department of Industrial Engineering, University of Tehran, Fooman, Rasht}
% \cormark[1]
\author[2]{Seyed Sajjad Fazeli}[]
% \fnmark[1]
\cormark[1]
\ead{sajjad.fazeli@wayne.edu}
% \ead[url]{www.cvr.cc, cvr@sayahna.org}

% \credit{Conceptualization of this study, Methodology, Software}

\address[2]{Department of Industrial and System Engineering, Wayne State university, Detroit, MI}

\author[3]{Soodabeh Milanlouei}[%
%   role=Co-ordinator,
%   suffix=Jr,
   ]
% \fnmark[2]
% \ead{cvr3@sayahna.org}
% \ead[URL]{www.sayahna.org}

% \credit{Data curation, Writing - Original draft preparation}

\ead{milanlouei.s@northeastern.edu}
\address[3]{Center for Complex Network Research, Northeastern University, Boston, MA}
% \ead[url]{www.cvr.cc, cvr@sayahna.org}

% \credit{Conceptualization of this study, Methodology, Software}
\author[2]{Ali Kamali Mohammadzadeh}[%
%   role=Co-ordinator,
%   suffix=Jr,
   ]
% \fnmark[2]
% \ead{cvr3@sayahna.org}
% \ead[URL]{www.sayahna.org}

% \credit{Data curation, Writing - Original draft preparation}

\ead{alikamali@wayne.edu}
% \address[3]{Center for Complex Network Research, Northeastern University, Boston, MA}
\author%
[4]
{Leyla Sadat Tavassoli}

\ead{Leylasadat.tavassoli@mavs.uta.edu}

\address[5]{Department of Industrial Manufacturing and Systems Engineering, University of Texas at Arlington, Arlington, TX}

% \author%
% [3]
% {Mirpouya Mirmozaffari}

% \ead{Mirpouya.mirmozaffari@mavs.uta.edu}

\cortext[cor1]{Corresponding author: sajjad.fazeli@wayne.edu}
% \cortext[cor2]{Principal corresponding author}
% \fntext[fn1]{This is the first author footnote. but is common to third
%   author as well.}
% \fntext[fn2]{Another author footnote, this is a very long footnote and
%   it should be a really long footnote. But this footnote is not yet
%   sufficiently long enough to make two lines of footnote text.}

% \nonumnote{This note has no numbers. In this work we demonstrate $a_b$
%   the formation Y\_1 of a new type of polariton on the interface
%   between a cuprous oxide slab and a polystyrene micro-sphere placed
%   on the slab.
%   }

\begin{abstract}
\abstract{Excessive greenhouse gas emissions from the transportation sector have led companies to move towards a sustainable supply chain network design. In this study we present a new bi-objective non-linear formulation where multiple inventory components are integrated into the location and routing decisions throughout the supply chain network. To efficiently solve the proposed model, we implement an exact method and four evolutionary algorithms for small and large-scale instances. Extensive computational results and sensitivity analysis are performed to validate the efficiency of the proposed approaches, both quantitatively and qualitatively. Besides, we run a statistical analysis to investigate whether there is any statistically significant difference between solution methods.}

\end{abstract}

% \begin{graphicalabstract}
% \includegraphics{figs/grabs.pdf}
% \end{graphicalabstract}

% \begin{highlights}
% \item Research highlights item 1
% \item Research highlights item 2
% \item Research highlights item 3
% \end{highlights}


\maketitle
\section{Introduction}
The intense competition in the global markets has forced companies to manage and design their supply chains in a more productive way. Designing a cost-efficient supply chain network asks for simultaneous decision making on strategic, tactical, and operation levels. Facilities' locations are considered as strategic decisions that strongly affect the tactical and operational aspects of a supply chain. Inventory policies and transportation activities are among the critical decisions at the tactical and operational levels. The integrity of these decisions is an essential factor that significantly reduces supply chain costs and leads to higher customer satisfaction. In recent decades, there have been extensive studies on the Location-Routing-Inventory problem (LRIP). The majority of these studies considered the same approach towards the inventory components. For example, the order quantity is assumed to be a constant across the supply chain network while, in reality, it could variate among distribution centers based on the customer's demand. The same situation could be observed for the Safety Stock (SS) level at different distribution centers. Besides, the inventory decisions could be significantly affected by uncertainty in distribution centers' demand. One of the first attempts to model a Location-Routing-Inventory problem was made by \cite{liu2003two}, where the authors considered a two-level supply chain consisting of customers and multiple depots with limited capacities. They incorporated the continuous review policy in their model. Later, the same model was solved by an innovative approach based on Tabu search and simulated annealing algorithms developed by \cite{liu2005heuristic}. Authors in \cite{shen2007incorporating} presented an innovative approach to incorporate non-linear inventory and routing costs into the facility location models considering randomness in customers' demand. The study by \cite{javid2010incorporating} introduced a novel mathematical model that simultaneously take into account the location, routing, and inventory decisions considering uncertainty in demand. The order quantity was the only inventory decision included in the model. \par Other than designing a cost-efficient supply chain network, due to the growing concern of global warming, companies strive to include the environmental aspects into their operations. It is reported that more than 28\% of the total greenhouse gas emissions come from the transportation sector in U.S., which is mainly released from road activities \cite{fazeli2020efficient,fazeli2020two}. Besides, distribution centers could be a source of carbon emission caused by power consumption and volume of inventory \cite{diabat2009carbon}.  Green logistics offers sustainable production and distribution strategies to the companies. Considering richer objectives and more operational constraints involved with sustainable logistics issues impose new challenges that lead to more complex combinatorial optimization problems. There are only a few studies that incorporated the environmental aspects into the LRIP. The research in \cite{zhalechian2016sustainable} developed a new multi-objective mathematical model that addresses routing, inventory, and location decisions. They validate their model by implementing a real case study and investigating various factors that affect $CO_{2}$ emissions in a supply chain network. \par LRIP is an extension of the facility location problem, and it is NP-Hard \cite{javid2010incorporating}. Numerous researchers proposed different methodologies to overcome the computational complexity of the LRIP. There are only a few studies that developed exact algorithms such as Benders Decomposition and Lagrangian relaxation to solve the LRIP \cite{yu2012large,zheng2019integrated}. The majority of literature focused on heuristic (\cite{javid2010incorporating,saragih2019heuristic,liu2003two,sajjadi2011multi}), meta-heuristic (\cite{yang2010research,nekooghadirli2014solving,forouzanfar2012using, wu2020hybrid}) and hybrid (\cite{ghorbani2016hybrid,zhalechian2016sustainable}) methods.
\par In this paper, considering uncertainties in Distribution Centers (DC) demand, lead time, and many other sources, we propose an innovative mathematical formulation that simultaneously considers location, inventory, and routing decisions under mixed uncertainty. The formulation specifically considers the important inventory decisions that arise in real-world applications which hasn't been fully investigated on the literature. Also, we take the green approach by adding a new objective function considering the vehicles' load and their emissions. To solve the proposed non-linear bi-objective mixed-integer model efficiently, we first linearize and then solve the small-scale instances with an exact method. For large-scale instances, we implement different meta-heuristic algorithms and statistically compare their performance to select the best solution approach .
\par The remainder of this paper is organized as follows:  Section
\ref{def} provides problem definition and mathematical formulation of the problem along with a subsequent reformulation. Section \ref{method} represents the solution methodologies, where we present exact and heuristic methods. Section \ref{results} provides details about the test instances, the numerical results, and statistical analysis for comparing the proposed approaches. Finally, Section \ref{con} provides concluding remarks. 

\section{Problem Definition and Model Formulation}\label{def}
\subsection{Problem Definition}
\begin{figure*}
\makebox[\textwidth]{
\includegraphics[scale=0.6]{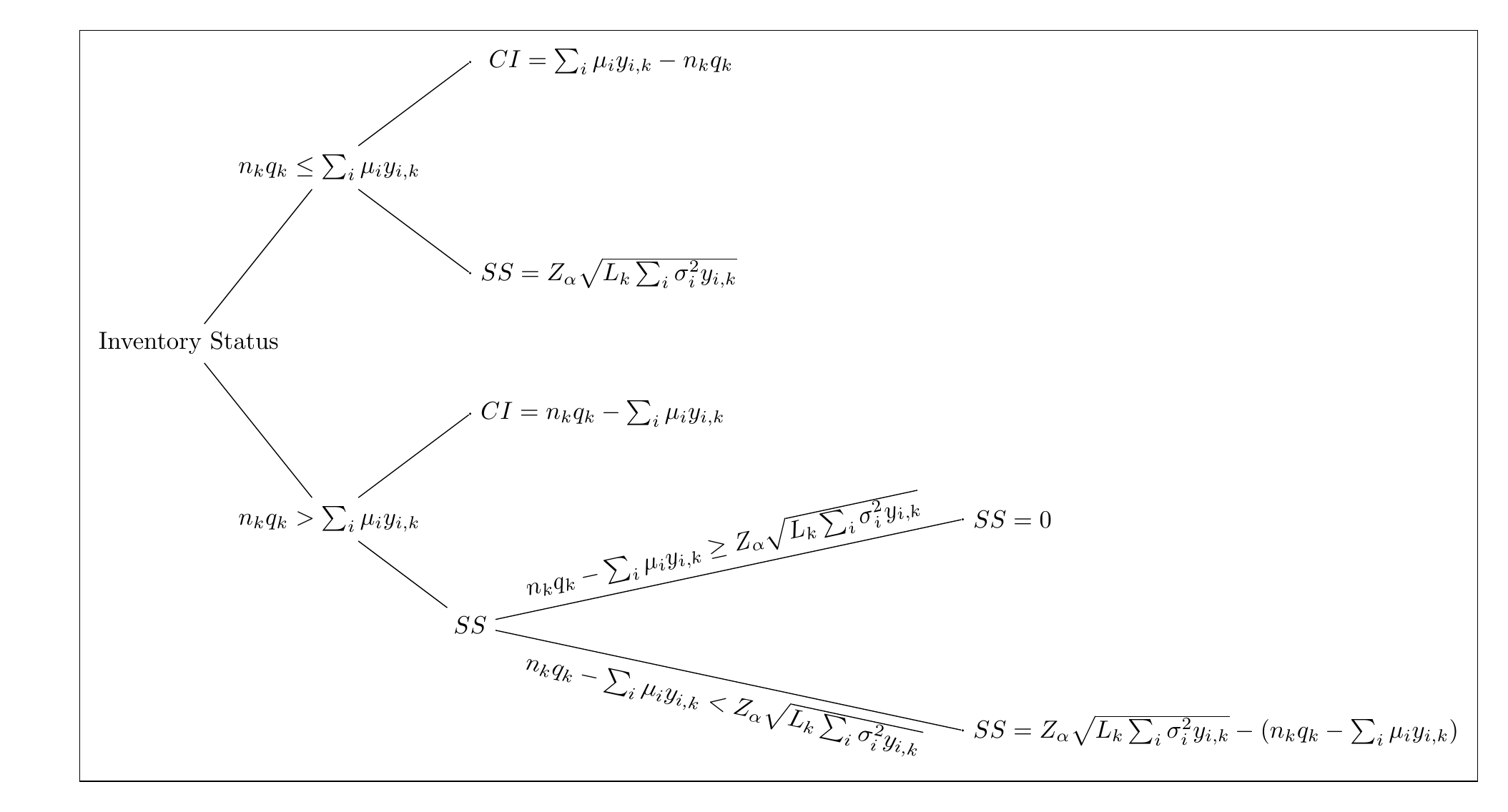}
}
\caption{Inventory status under different scenarios}
\label{Inv_Scenarios}
\end{figure*}

We consider a three-echelon supply chain network where the first echelon contains a supplier that receives the orders from the distribution centers. The second level includes distribution centers, where they receive orders from retailers in the third echelon. Once the distribution centers receive orders from the retailers, they send the order to the supplier based on the continuous inventory policy. The continuous policy for DC $k$ is denoted as $(q_{k}, r)$. In this policy, the inventory of DCs is continuously reviewed. Whenever the inventory amount becomes less than or equal to the reorder point ($r$),  an order of size $q_k$ is placed to the supplier. In reality, the demand in a supply chain is affected by many sources of uncertainty. Authors in \cite{montgomery2009engineering} suggested that the normality assumption is a reasonable approximation for sufficiently large demand values.  We assume that the demand of retailer $i$ follows normal distribution $D_{i}  \sim N(\mu_{i},\sigma_{i})$. We define $n_k$ as order frequency and $q_k$ as the order size of DC $k$. So, the annual inventory received by  distribution $k$ is $n_k q_{k}$. Considering multiple retailers, the annual average demand and variance for all retailers assigned to DC $k$ are equal to $\sum_{i}\mu_{i}y_{i,k}$ respectively and $\sum_{i}\sigma^{2}_{i}y_{i,k}$ where $y_{i,k}$ is  a binary variable equal to one, if retailer $i$ is assigned to DC $k$.
Due to uncertainty in demand, the value of $n_{k} q_{k}$ and $\sum_{i}\mu_{i} y_{i,k}$ could be dissimilar which will result in different scenarios for the SS and Cycle Inventory (CI). Figure \ref{Inv_Scenarios} shows the inventory status under different scenarios. The $L_k$ denotes the lead time for distribution $k$.\\
The goal is to locate a set of distribution centers, determine the inventory decisions, and select the best possible routes for a fleet of vehicles to satisfy the retailers demand such that the total cost of the system is minimized.
The main assumptions used to formulate the model are as follow:
\begin{itemize}
\item Each retailer has an uncertain demand which follows a normal distribution.
\item Each distribution center follows continuous review policy $(q_k, r_{k})$.
\item Vehicles are different in terms of capacity (heterogeneous fleet).
\item Vehicles return to the start point after delivering the products.
\item The inbound transportation (from supplier to distribution center) is considered to be direct (no routing), and the outbound transportation (from distribution to retailers) is decided based on routing strategy. 
\item One type of product is distributed in the network. 
\item The amount of $CO_2$ emission depends on the traveled distance, the load weight of vehicles, and the DC emissions caused by volume of flow.
\item There is only one supplier in the network denoted as $s_{0}$.
\end{itemize}
\subsection{Notation}
\begin{itemize}
	\item Sets
	\end{itemize}
	\begin{itemize}[label=-]
		\item $K$: Set of distribution centers, with $k\in K$
		\item $I$: Set of retailers, with $i \in I$
		\item $V^{in}$: Set of vehicles used in the inbound transportation, with $v \in V^{in}$
		\item $V^{out}$: Set of vehicles used in the outbound transport ion (transport from suppliers to DCs), with $v \in V^{out}$
		\item $V$: Set of all vehicles, with $v \in V$ and $V= V^{in}\cup V^{out}$.
		\item $S$: Set of all nodes including retailers and distribution centers with $s \in S=K \cup I $
	\end{itemize} 
	\begin{itemize}
	\item Model parameters
	\end{itemize}
	\begin{itemize}[label=-]
		\item $f_{k}$: Fixed establishment cost for DC $k \in K$
		\item $g_{k}$: Fixed transportation cost from supplier to DC $k \in K$.
		\item $a_{k}$: Fixed ordering cost of DC $k$ to the supplier
		\item $b_{k}$: Supply cost from supplier $s_0$ to DC $k$
		\item $d_{s_{0},k}$: Distance between supplier  $s_0$ and DC $k \in K$
		\item $l_{k}$ :Lead time of DC $k \in K$
		\item $\mu_{i}$: Mean of annual demand for retailer  $i \in I$
		\item $\sigma_{i}^2$: Variance of annual demand for retailer  $i \in I$
		\item $d_{i,j}^{'}$: Distance between node $i$ and node $j$, with $i,j \in S$
		\item $c_{i,j}$: Transportation cost between node $i$ and node $j$, with $i,j \in S$
		\item $h_k$: Inventory holding cost of DC $k$ with $k \in K$
		\item $\rho_{v}^{0}$: Amount of fuel consumption by vehicle $v \in V$ with no load
		\item $\rho_{v}^{*}$: Amount of fuel consumption by vehicle $v \in V$ with maximum load capacity
		\item $o_{v}$: Capacity of vehicle $v \in V$
		\item $o'_{k}$: Annual capacity of DC $k \in K$
		\item $f_{v}$: Fixed cost for each vehicle $v \in V$
		\item $\gamma_{v}$: Carbon emission of vehicle $v \in V$
		\item $\eta_k$: Weight factor associated with carbon emission from DC $k$
		\item $\alpha$: Service level of distribution centers
		\item $\theta$: Weight factor associated with inventory cost
		\item $\beta$: Weight factor associated with shipping cost
% 		\item $\eta$: \textcolor{blue}{What is this?} weight factor associated with carbon emission.
		\item $M$: Sufficiently large number
\end{itemize}
\begin{itemize}
	\item Decision variables
	\end{itemize}
	\begin{itemize}[label=-]
	\setlength\itemsep{1em}
		\item $x_{k}$: 1 if DC $k \in K$ is selected for establishment; 0 otherwise.
		\item $q_{k}$: Quantity of orders by DC $k \in K$
		\item $n_{k}$: Number of orders by DC $k \in K$ during a year
		\item $y_{i,k}$: 1 if retailer $i \in I$ is selected to be sourced by DC $k \in K$; 0 otherwise.
		\item $r_{k,v}$: 1 if vehicle $v \in V$ is assigned to DC $k \in K$ for transportation; 0 otherwise.
		\item $w_{i,j,v}$: 1 if edge (i,j) is visited by vehicle $v \in V$, with $i,j \in S$; 0 otherwise.
		\item $u_{i,j,v}$: Flow of products from node $i$ to node $j$ carried by vehicle $v$, with $i,j \in S$
		\item $m_{i,v}$: Auxiliary binary variable for removing sub-tours related to retailer $i \in I$ and vehicle $v \in V$
		\item $t_{k}$: Auxiliary binary variables associated with DC $k \in K$
		\item $t_{k}^{'}$: Auxiliary binary variables associated with DC $k \in K$
	\end{itemize}
\subsection{Mathematical Formulation}
The Non-linear bi-objective mixed integer linear programming model is defined as follows:
% \vspace{-10cm}
\begin{alignat}{2}& 
\mathrm{ Min } \nonumber\\
%==================================================
%           First Objective Function
%==================================================
&\boldsymbol{z_{1}} = 
\underbrace{\sum_{k \in K}\Bigg(f_{k} x_{k} +\sum_{v \in V^{in}}f_{v}w_{s_{0},k,v} +\sum_{i \in I}\sum_{v \in V^{out}}f_{v}w_{k,i,v}\Bigg)}_{\text{Fixed DC establishment and fleet cost}}+ \nonumber\\ 
&\underbrace{\beta \cdot \sum_{k \in K}\Bigg(\sum_{v \in V^{in}} (g_{k}+b_{k}q_{k})n_{k} + \sum_{v \in V^{out}}\sum_{i,j \in S}c_{i,j}w_{i,j,v}r_{k,v}n_{k}\Bigg)}_{\text{Total shipping cost}} + \nonumber\\
&\underbrace{\sum_{k \in K} a_{k}n_{k}}_{\text{Total fixed ordering cost}}+ \underbrace{\theta \cdot \sum_{k \in K}h_{k}\Bigg(t_{k}(\sum_{i \in I}\mu_{i}y_{i,k}-n_{k}q_{k})}+\nonumber\\
&\underbrace{t_{k}Z_{\alpha} \sqrt{l_{k}\sum_{i \in I}\sigma_{i}^{2}y_{i,k}}
+(1-t_{k})(n_{k}q_{k}-\sum_{i \in I}\mu_{i}y_{i,k}) +}\nonumber\\ &\underbrace{+(1-t_{k})(1-t_{k}^{'})\big(Z_{\alpha} \sqrt{l_{k}\sum_{i \in I}\sigma_{i}^{2}y_{i,k}}-(n_{k}q_{k}-\sum_{i \in I}\mu_{i}y_{i,k})\big)\Bigg)}_{\text{Total inventory cost}}&\label{obj1}
%==================================================
%           Second Objective Function
%==================================================
\end{alignat}
\begin{alignat}{2}
& \boldsymbol{z_{2}} = \underbrace{\sum_{v \in V^{in}}\sum_{k \in K}\gamma_{v}\big(2\rho_{v}^{0}w_{s_{0},k,v} + (\frac{\rho_{v}^{*}-\rho_{v}^{0}}{o_{v}})u_{s_{0},k,v}\big)n_{k}d_{s_{0},k}}_{\text{Total $CO_{2}$ emission in inbound transportation}} + \nonumber\\
& \underbrace{\sum_{v \in V^{out}}\sum_{k \in K}\sum_{i,j \in S} \gamma_{v} \big(\rho_{v}^{0}w_{i,j,v}+(\frac{\rho_{v}^{*}-\rho_{v}^{0}}{o_{v}})u_{i,j,v}\big)n_{k}d_{i,j}^{'}r_{k,v}}_{\text{Total $CO_{2}$ emission in outbound transportation}}+\nonumber\\
&\underbrace{\eta_k \cdot \sum_{k \in K}h_{k}\Bigg(t_{k}(\sum_{i \in I}\mu_{i}y_{i,k}-n_{k}q_{k})+ t_{k}Z_{\alpha} \sqrt{l_{k}\sum_{i \in I}\sigma_{i}^{2}y_{i,k}}
}+ \nonumber\\
&\underbrace{(1-t_{k})(n_{k}q_{k}-\sum_{i \in I}\mu_{i}y_{i,k}) +}\nonumber\\ &\underbrace{+(1-t_{k})(1-t_{k}^{'})\big(Z_{\alpha} \sqrt{l_{k}\sum_{i \in I}\sigma_{i}^{2}y_{i,k}}-(n_{k}q_{k}-\sum_{i \in I}\mu_{i}y_{i,k})\big)\Bigg)}_{\text{Total $CO_2$ emissions from DCs}}\label{obj2}\\
& \text{s.t. }  \nonumber\\
%==================================================
%           Constraints
%==================================================
% Constraint3
& x_{k}\geq y_{i,k} \quad  \forall k \in K, i \in I  \label{cons3} \\
% Constraint4
& \sum_{k\in k}y_{i,k} = 1 \quad \forall i \in I  \label{cons4} \\
% Constraint5
& \sum_{j \in S}\sum_{v \in V^{out}}w_{i,j,v}=1 \quad \forall i \in I  \label{cons5} \\
% Constraint6
& \sum_{v \in V^{in}}w_{s_{0},k,v} \geq x_{k} \quad \forall k \in K  \label{cons6} \\
% Constraint7
& m_{i,v}-m_{j,v}+(|I|\cdot w_{i,j,v})\leq |I|-1 \quad \forall i,j \in I,v\in V^{out}\label{cons7} \\
% Constraint8
& \sum_{j \in S}w_{i,j,v}-\sum_{j \in S}w_{j,i,v}=0  \qquad \forall i \in S,v \in V^{out}  \label{cons8} \\
% Constraint9
& \sum_{k \in K}w_{k,s_{0},v}\leq 1 \qquad \forall v \in V^{in} \label{cons9}\\
% Constraint10
& \sum_{k \in K}\sum_{i \in I}w_{k,i,v}\leq 1 \qquad \forall v \in V^{out}  \label{cons10}\\ 
% Constraint11
& \sum_{j \in S}w_{i,j,v}+\sum_{j \in S}w_{k,j,v}-y_{i,k}\leq 1 \qquad \forall k \in K, i \in I,v \in V^{out} \label{cons11} \\
&  \sum_{j \in S}\sum_{v \in V^{out}} \big(u_{i,j,v} - u_{j,i,v}\big)\leq M (1- y_{i,k})+\frac{\mu_{i}}{n_{k}}  \quad \forall i \in I,k \in K \label{cons12_1} \\
&  \sum_{j \in S}\sum_{v \in V^{out}} \big(u_{i,j,v} - u_{j,i,v}\big)\geq -M (1- y_{i,k})+\frac{\mu_{i}}{n_{k}}  \quad \forall i \in I,k \in K \label{cons12_2} \\
% Constrain13
& u_{s_{0},k,v}\leq o_{v} w_{s_{0},k,v} \quad \forall k \in K,v \in V^{in} \label{cons13} \\
% Constraint14
& u_{i,j,v}\leq o_{v} w_{i,j,v} \quad  \forall i,j \in S,v \in V^{out} \label{cons14}\\
% Constraint15
& \sum_{i \in I}\mu_{i} \cdot y_{i,k}  \leq o_{k}^{'} \quad  \forall k\in K\label{cons15}\\
% Constraint16
& q_{k}=\sum_{v \in V^{in}} u_{s_{0},k,v} \quad  \forall k\in K \label{cons16}\\
% Constraint17
& r_{k,v}=w_{s_{0},k,v} \quad  \forall k\in K,v \in V^{in} \label{cons17}\\
% Constraint18
& r_{k,v}=\sum_{i \in I}w_{k,i,v} \quad  \forall k\in K,v \in V^{out} \label{cons18}\\
% Constraint19
& Mt_{k}\geq  \sum_{i \in I}\mu_{i}y_{i,k} - n_{k}q_{k} \quad  \forall k\in K \label{cons19}\\
% Constraint20
&-M(1-t_{k})\leq  \sum_{i \in I}\mu_{i}y_{i,k} - n_{k}q_{k}  \quad  \forall k\in K  \label{cons20}\\
% Constraint21
&Mt_{k}^{'}\geq n_{k}q_{k}-\sum_{i \in I}\mu_{i}y_{i,k}-Z_{\alpha}\sqrt{l_{k}\sum_{i \in I}\sigma_{i}^{2}y_{i,k}} \quad \forall k\in K  \label{cons21}\\
% Constraint22
&-M(1-t_{k}^{'})\leq n_{k} q_{k} - \sum_{i \in I}\mu_{i} y_{i,k}- Z_{\alpha} \sqrt{l_{k}\sum_{i \in I}\sigma_{i}^{2} y_{i,k}} \quad \forall k\in K \label{cons22} \\
% Constraint23
& x_{k} \in \mathbb{Z}, y_{i,k},w_{i,j,v},t_{k},t_{k}^{'}, r_{k,v} \in \{0,1\} \forall k \in K,i,j \in S,v \in V  \label{cons23}\\
% Constraint24
& n_{k},q_{k},u_{i,j,v},m_{i,v} \geq 0 \quad \forall k \in K,i,j \in S,v \in V  \label{cons24}
\end{alignat}
In equation (\ref{obj1}), the first objective function minimizes the total costs of the network, including fixed annual DCs' establishment cost, fixed ordering cost, and fixed vehicle's cost and expected inventory cost. The second objective function (\ref{obj2}) minimizes the total carbon emission from inbound and outbound transportation activities, also, the carbon emissions from the distribution centers caused by flow of inventory.
Constraints (\ref{cons3}) ensure that the retailers are only assigned to the established DCs. Constraints (\ref{cons4}) guarantee that each retailer is sourced by only one DC. Constraints (\ref{cons5}) ensure that the demand of each customer is delivered by only one vehicle. Constraints (\ref{cons6}) ensure that the vehicles dispatched from supplier to DC if only DC is established. Constraints (\ref{cons7}) and (\ref{cons8}) remove sub-tours and maintain flow in outbound transportation. Constraints (\ref{cons9}) and (\ref{cons10}) make sure that at most one established DC is included in each tour. Constraints (\ref{cons11}) preserve the connection between the location and routing decisions, meaning that if vehicle $v$ is dispatched from DC $k$ and visits retailer $i$, then DC $k$ should be assigned to retailer $i$. Constraints (\ref{cons12_1}) and (\ref{cons12_2}) balance the product flow and the demand from DCs to retailers. Constraints (\ref{cons13}) and (\ref{cons14}) impose the vehicles' capacity limitations in inbound and outbound transportation. Constraints (\ref{cons15}) ensure the DC's capacity is not violated. Constraints (\ref{cons16}) guarantee that the total shipment by vehicles in inbound transportation is equal to the order quantity. Constraints (\ref{cons17}) and  (\ref{cons18}) state that the vehicle $v$ can be dispatched if it is assigned to DC $k$. Constraints (\ref{cons19}) to (\ref{cons22}) are the mathematical representation of Figure \ref{Inv_Scenarios}. Constraints (\ref{cons23}) to (\ref{cons24}) are variables restrictions.
%%%%%%%%%%%%%%%%%%%%%%%%%%%%%%%%%%%%%%%%%%
\section{Methodology}\label{method}
The model proposed in the previous section is a non-linear bi-objective mixed-integer formulation. To efficiently solve the problem,  we first linearize the formulation and then offer an exact method to solve the small-scale instances. For large scale instances, we suggest four different evolutionary algorithms.
\subsection{Exact Method}
One of the approaches that is often used to obtain the Pareto solutions in multi-objective problems is augmented $\varepsilon$-constraint (A$\varepsilon$-c) method \cite{cohon2004multiobjective}. In the $A\varepsilon$-constrain method, one of the objectives is selected as the main objective and the others ones will be added to the constraint sets. Suppose that there are $m$ objective functions, then we have:
\begin{alignat}{3}
&\text{min} f_{1}(x) \nonumber\\
&\text{s.t}\nonumber\\
& f_{2}(x) \leq e_2,...,f_{m}(x)\leq e_m \\
& x \in S
\end{alignat}
where $x$ is a vector of decision variables and $S$ is the feasible region. We can obtain the set of solutions by changing the right-hand side of the constrained objective functions ($e_i$). However, the drawback of the $\varepsilon$ constraint is that once the optimal solution for the main objective is obtained, the method has no obligation to find the best value for other objectives. In order to resolve this issue, \cite{mavrotas2013improved} suggested using a lexicographic optimization technique to ensure multiple Pareto solutions. The interested readers may refer to \cite{mavrotas2013improved} for more details. To efficiently implement the (A$\varepsilon$-c), we reformulate the proposed model to linear mixed integer problem. The details about the linearization methods are provided in  \ref{linearization}.
\subsection{Solution Representation}
\begin{figure}
		\centering
		\includegraphics[scale=0.3]{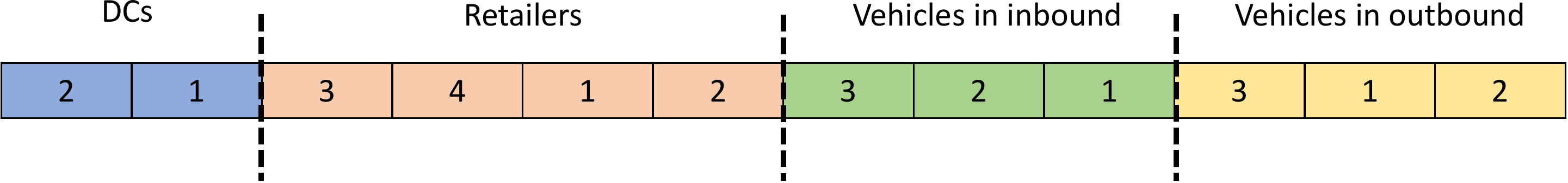}
		\captionsetup{justification=centering}
		\caption{A sample chromosome ordered based on the initial values}
	\label{S2}
	\end{figure}

The first step in implementing a meta-heuristic algorithm is the solution representation. The designed chromosome for this problem is a  string of random numbers between $(0,1)$ with the length of $|K|+|I|+|V^{in}|+|V^{out}|$. For example,  a sample chromosome with two DCs, four retailers, and three vehicles on both outbound and inbound routes is as follows:
\begin{figure}
\centering
		\includegraphics[scale=0.3]{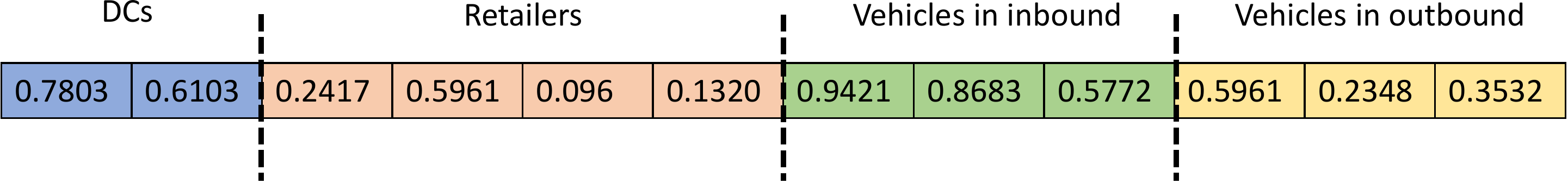}
		\captionsetup{justification=centering}
		\caption{Sample chromosome designed for solution representation}
		\label{S1}
	\end{figure}
In the next step, sub-strings are ordered based on their values to determine each cell's priority. For example, in Figure \ref{S2}, the first sub-string indicates the priority of each DC for establishment. The second sub-string illustrates the priority of the assignment of retailers to the DCs. The third and fourth sub-string prioritize the allocation of vehicles.
To obtain a feasible solution, we designed an efficient decoding strategy. According to the chromosome in Figure \ref{S2}, we select DC 2 for establishment. We assign the retailers DC 2 based on their order   until DC’s capacity is violated. Then, we choose the next DC and continue as before until all the retailers are assigned to the DCs. We make sure that the assignments preserve DCs capacity constraints. Let’s assume that the capacity of  DC 1 is 1000 units and the capacity of retailers 1 to 4  are 300, 700, 400, and 600 units, respectively. Thus, retailers 3 and 4 (based on their order in the chromosome) are assigned to the DC 2 and retailers 1 and 2 to DC 1.It should be noted that if the retailers order were [3,1,2,4] instead of [3,4,1,2], retailers 3 and 1 would have been assigned to DC 2, and retailer 1 to DC 1. Due to capacity violation, we could not have assigned retailer 4 to any of the DCs. So, this solution is considered to be an infeasible solution. We penalize the infeasible solutions with cost with respect to the number of unassigned retailers. Once the retailers' assignment to DCs is determined, the mean and variance of demand of each DC are calculated based on the aggregated demand of retailers.The next step is the assignment of inbound and outbound transportation to the DCs. The number of vehicles is dependent on the order quantity and number of orders by each DC. Following the above example, the annual demand for DC 2 is 400+600=1000 units. If we consider one-time order during a year, the entire order should be delivered to the DC and then to retailers in one shipment. Suppose that the capacity of vehicles in inbound transportation is 300 units. So, it is impossible to deliver the products to the DC 2 with three vehicles. We denote $n^{max}_{k}$ as the maximum number of orders by DC $k$.  So, for the different number of orders ($n$), we can estimate the sufficient number of vehicles. For example, for $n=2$ for DC 2, the order quantity is 500, and the minimum number of vehicles is two. The method is repeated for other DCs to obtain the minimum number of vehicles. Based on the above example, the estimation of number of vehicles needed for each DC is provided in Table \ref{example_0}. 
\begin{table}
\centering
\begin{tabular}{|c|cccccc|}
\hline
     & $n=1$ & $n=2$ & $n=3$ & $n=4$ & $n=5$ & $n=6$ \\ \hline
DC 1 & -   & 2   & 2   & 1   & 1   & 1   \\ 
DC 2 & -   & 2   & 2   & 1   & 1   & 1   \\ \hline
\end{tabular}
\caption{Estimation for the number of vehicles at each DC based on different number of orders}
\label{example_0}
\end{table}
We implement the same method for the outbound vehicles to create a similar table to Table \ref{example_0}. Once the tables for both inbound and outbound transportation are formed, for every DC, we randomly select an $n$ such that it is feasible in both layers of transportation. 
In the next step, we assign the vehicles to the distribution based on their order in the chromosome. Suppose that, $n_{k=2}=2$ and  $n_{k=1}=4$, then vehicles 3 and 2 are assigned to DC 2 and vehicle 1 to DC 1. Since the amount of orders for DC 2 in each shipment is 500 units, the amount of load for vehicle 3 and 2 are 300 and 200 units, respectively. Also, vehicle 1 is assigned to DC 1 with 125 units of load in each shipment. 
Once all decisions are made, the amount of $CO_{2}$ is calculated based on the vehicles’ load.
\subsection{Evolutionary Algorithms}
Evolutionary algorithms are stochastic techniques inspired by the process of natural selection, which are extensively applied to solve different classes of NP-hard problems.
In the literature, there are several evolutionary algorithms developed for multi-objective optimization problems. These algorithms use the Pareto dominance concept to lead the search process and return the best solutions in the form of a Pareto optimal set, while preserving the convergence and the diversity in the solution set.
To solve the proposed model in Section \ref{def} efficiently, we implement four meta-heuristic algorithms with different characteristics. These algorithms are Non-dominated Sorting Genetic Algorithms II (NSGA-II), Non-dominated Ranked Genetic algorithm (NRGA), Strength Pareto Evolutionary Algorithm II(SPEA-II), and Pareto Envelope-based Selection Algorithm II(PESA-II). We intend to evaluate the performance of well-known population-based algorithms for the proposed multi-objective location-inventory-routing problem.
\subsubsection{Non-dominated Sorting Genetic Algorithm II (NSGA-II)}
Proposed by \cite{deb2002fast} as a biological heuristics algorithm, NSGA-II is among the most widely used multi-objective evolutionary algorithms. The algorithm utilizes a population of individuals and employs an elitism based sorting method. NSGA-II consists of two main operators, namely non-dominated sorting and crowding distance procedures. The non-dominated sorting method ranks all the solutions into different non-dominated levels according to the Pareto dominance principle. The crowding distance method, on the other hand, preserves the diversity of the solutions by calculating the dispersion of any two neighboring solutions in each front. These two procedures shape the Pareto front at each iteration \cite{tavana2018evolutionary,tavassoli2020integrated}.
After generating an initial parent population ($P_0$), all non-dominated individuals are sorted. Using the tournament selection strategy, crossover, and mutation operators, the offspring population ($Q_0$) is generated afterward. Combining the parent and offspring populations, the combined population ($R_t= P_t \cup Q_t$ ) is formed at each generation $t$. The Pareto front is obtained by non-dominated sorting of combined population and estimation of the crowded distance of solutions. The subsequent parent population ($P_{t}+1$) is created by selecting the best individuals according to the rank and crowded distance. This procedure continues until the termination criterion is met.
\subsubsection{Non-dominated Ranked Genetic Algorithm (NRGA)}
Introduced by \cite{al2008non}, NRGA is analogous to NSGA-II with one exception in the selection mechanism. Instead of the tournament selection operator utilized in NSGA-II, NRGA employs a ranked-based roulette wheel (RBRW) selection operator, combined with a Pareto-based population-ranking algorithm. According to RBRW, a parent i is selected with the following probability:
$$ P_{i}=\frac{2\times Rank_{i}}{N\times (N+1)}$$
Where $P_i$ and $N$ indicate the probability of being selected and the number of individuals in the population, respectively.
\subsubsection{Strength  Pareto Evolutionary Algorithm II (SPEA-II)}
Presented in 1999 and 2001, SPEA and SPEA-II are considered as extensions of the genetic algorithm. SPEA-II employs three main operators: a fine-grained fitness assignment strategy, a nearest neighbor density estimation method, and an enhanced archive truncation method that guarantees border solutions are preserved \cite{lara2019multiobjective}
This fitness assignment procedure has two main parts: the raw fitness based on the concept of dominance and a density estimation based on the k-nearest neighbor method. Besides, a truncation method is utilized to keep a fixed number of individuals in the archive. The interested reader may refer to \cite{zitzler2001spea2} for more detailed descriptions.
\subsubsection{Pareto Envelope-based Selection Algorithm (PESA-II)}
PESA-II is an evolutionary optimization algorithm for solving multi-objective optimization problems developed by \cite{corne2001pesa}.
This algorithm is a modified version of its predecessor, PESA, in which region-based selection is used for assigning selective fitness. The algorithm maintains two populations: a fixed size internal population, and an external population, also called the archive set. Containing only the non-dominated solutions, the archive is being maintained for the selection at each iteration. The core of the algorithm can be summarized in six main steps as follows:
\begin{enumerate}
    \item initialize the exterarchive and evaluate the generated internal population
    \item incorporate the initial population (non-dominated members from the archive) into the archive.
    \item If termination criteria are met, then stop. Otherwise, erase the internal population and reiterate step 4 until the generation of the new solution.
    \item From the archive, select new parents and produce a new child using crossover and mutation. Set the probability of this process as $p$
    \item Make the selection of one parent with the probability $1-p$, perform the mutation of this parent to produce a child.
    \item Return to step 2.
\end{enumerate}
Now a candidate solution may enter the archive if and only if it has non-dominancy over the internal population. Any dominated candidate is removed from the archive.
\section{Experiments and Results}\label{results}
In this section, we compare the computational performance of the exact method with the four meta-heuristic algorithms explained in section \ref{method}. Experiments were implemented in GAMS software equipped with the CPLEX 12.8 solver using a laptop with Intel(R) Core(TM) i7-9750H CPU @ 2.60 GHz, and 32GB RAM. It is worth to mention that the values used in the model are generated from probability distribution provided in Table \ref{params}.
\begin{table*}
	\centering
	\begin{tabular}{|c|  c  c c c c c c c |} 
		\hline
		Parameters &$f_{k}$&$g_{k}$&$a_{k}$&$b_{k}$& $L_{k}$&$h$& $\mu_{i}$&	$\sigma_{i}^2$ \\
		\hline
		 Value& U(500-1000)&U(10-15)&U(10-15)& U(5-10)& U(6-10)& U(5-10)&U(400-1500) & U(10-100)\\
		\hline
	\end{tabular}
	\caption{Model's Parameters}
	\label{params}
\end{table*}
\subsection{Parameter Tuning}
The performance of the meta-heuristic algorithms is highly dependent on the input parameters. One of the most efficient methods for tuning the algorithm parameters is the Taguchi method \cite{taguchi1986introduction}. The Taguchi can return a large amount of information with the least number of experiments. The Taguchi method makes use of orthogonal arrays, which approximate the effects of factors on the response mean and variation. Factors such as Noise factors result in variability in the performance of a system and cannot be controlled during production. The goal of the Taguchi method is to reduce the effect of uncontrollable factors and determine the best level(with higher signal-to-noise ($S/N$) ratio) for controllable factors. The $S/N$ is calculated as follows:
$$S/N= -10\cdot log\bigg(\frac{1}{n}\sum_{i=1}^{n}y^{2}_{i}\bigg)$$
Since in multi-objective problems, a set of non-dominated solutions is considered optimal solutions, various metrics are often used to compare the solution algorithms. In this study, we use the following metrics:
\begin{itemize}
    \item \textbf{Quality Metric (QM)}: Once non-dominated solutions from the four algorithms are obtained and stored in an archive, each pair of solutions are compared to each other and dominated solutions are eliminated. The algorithm's share in the archive shows its quality.
    \item \textbf{Spacing Metric (SM)}: The SM metric provides valuable information regarding the distribution of non-dominated solutions in the solution space. The SM is calculated as follows:
    $$SM = \frac{\sum^{n-1}_{i=1}|\Bar{d}-d_{i}|}{(n-1)\Bar{d}}$$
    where $n$ is the number of Pareto solution, $d_i$ is the euclidean distance between successive Pareto solutions, and $\Bar{d}$ is the average of $d_i$s. An algorithm with lower SM is more preferable.
    \item \textbf{Mean Ideal Distance(MID)}: The MID represents the distance between the best solution and Pareto solutions and calculated as follows:
    $$MID=\frac{\sum_{i=1}^{n}\sqrt{\big(\frac{f_{1,i}-f_{1}^{best}}{f_{1,total}^{max}-f_{1,total}^{min}}\big)^{2}+\big(\frac{f_{2,i}-f_{2}^{best}}{f_{2,total}^{max}-f_{2,total}^{min}}\big)^{2}}}{n}$$
    where $f_{i,total}^{max}$ and $f_{i,total}^{min}$ are maximum and minimum values among the Pareto solutions for objective $i$ respectively. The $f_{i}^{best}$ denotes the best solution for objective $i$.
    \item \textbf{Diversification Matrix (DM)}: The DM shows the diversity in the Pareto solutions which is calculated as follows:
    $$DM=\sqrt{\sum_{i=1}^{n} max(|x_{i}-y_{i}|) }$$
    where $x_i$ and $y_i$ denote the Pareto solutions for objective $i$.
\end{itemize}
To determine efficient values for algorithm parameters, we implemented a three-level Taguchi design in Minitab 17, including 9 experiments for NSGA-II and NRGA and 27 experiments for SPEA-II and PESA-II input parameters. Table \ref{tuned} shows the tuned parameters. The details regarding the Taguchi experiments is provided in \ref{Taguhci}.
\begin{table} 
	\caption{\hspace{2.2cm} Tuned Parameters for each algorithm}
	\centering
	\begin{tabular}{|c| c c c c|} 
		\hline
		\multirow{1}{*}{Parameters}& NSGA-II & NRGA & SPEA-II & PESA-II\\
		\hline
		\multirow{1}{*}{Population Size} & 100 & 150 & 100 &100  \\
  	    \multirow{1}{*}{Archive Size} & - & - & 300 &100 \\
  	    \multirow{1}{*}{Crossover Percentage} &0.7	&	0.7	&	0.9	&	0.7\\
  	    \multirow{1}{*}{Mutation Percentage} &0.3	&	0.2	&	0.2	&	0.2\\
  	    \multirow{1}{*}{Mutation Rate} &0.03	&	0.05	&	0.03	&	0.05\\
  	    \multirow{1}{*}{Selection Pressure} & -	&	-	&	-	&	3\\
  	    \multirow{1}{*}{Deletion  Pressure} & -	&	-	&	-	&	3\\
  	    \hline
	\end{tabular}
	\label{tuned}
\end{table}
\subsection{Instance Generation}
To demonstrate algorithms' performance, we generated 12 test problems where the first five are considered small-scale and the rest are large-scale problems. Table \ref{instances} provides the details about the instances, where $|k|$,$|I|$,$|V^{in}|$, and $|V^{out}|$ indicate the number of DCs, retailers, vehicles in inbound, and vehicle in outbound transportation, respectively.
\subsection{Results and Discussion}
\begin{figure}
\centering
\includegraphics[scale=0.7]{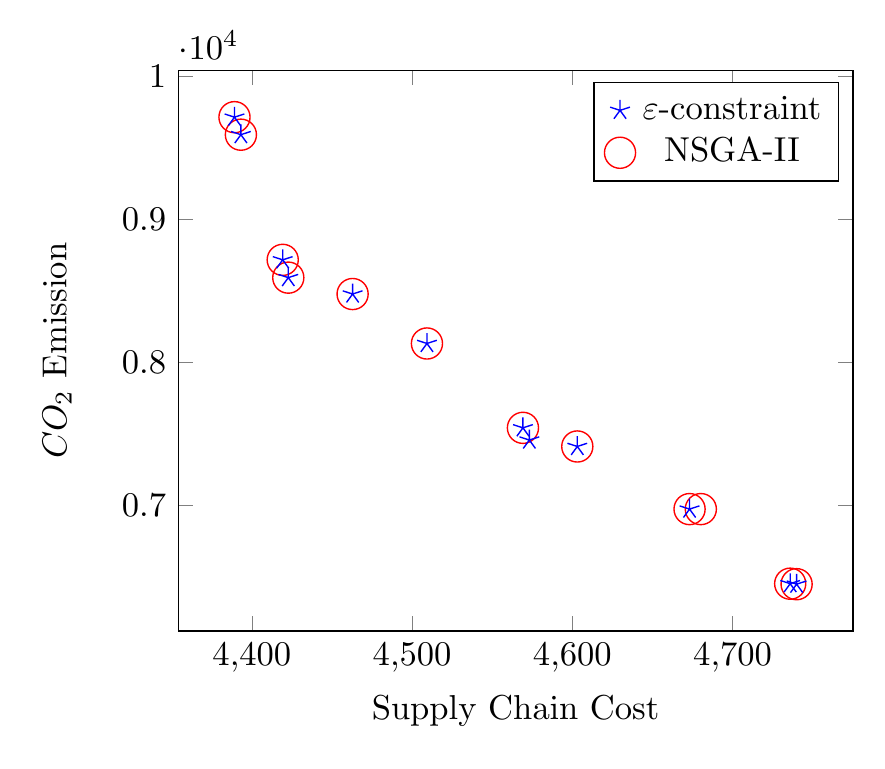}
\captionsetup{justification=centering}
\caption{Comparison between solutions found by NSGA-II algorithm and $\varepsilon$-c for test problem 5}
\label{pareto_results}
\end{figure}

\begin{figure}
\centering
\includegraphics[scale=0.6]{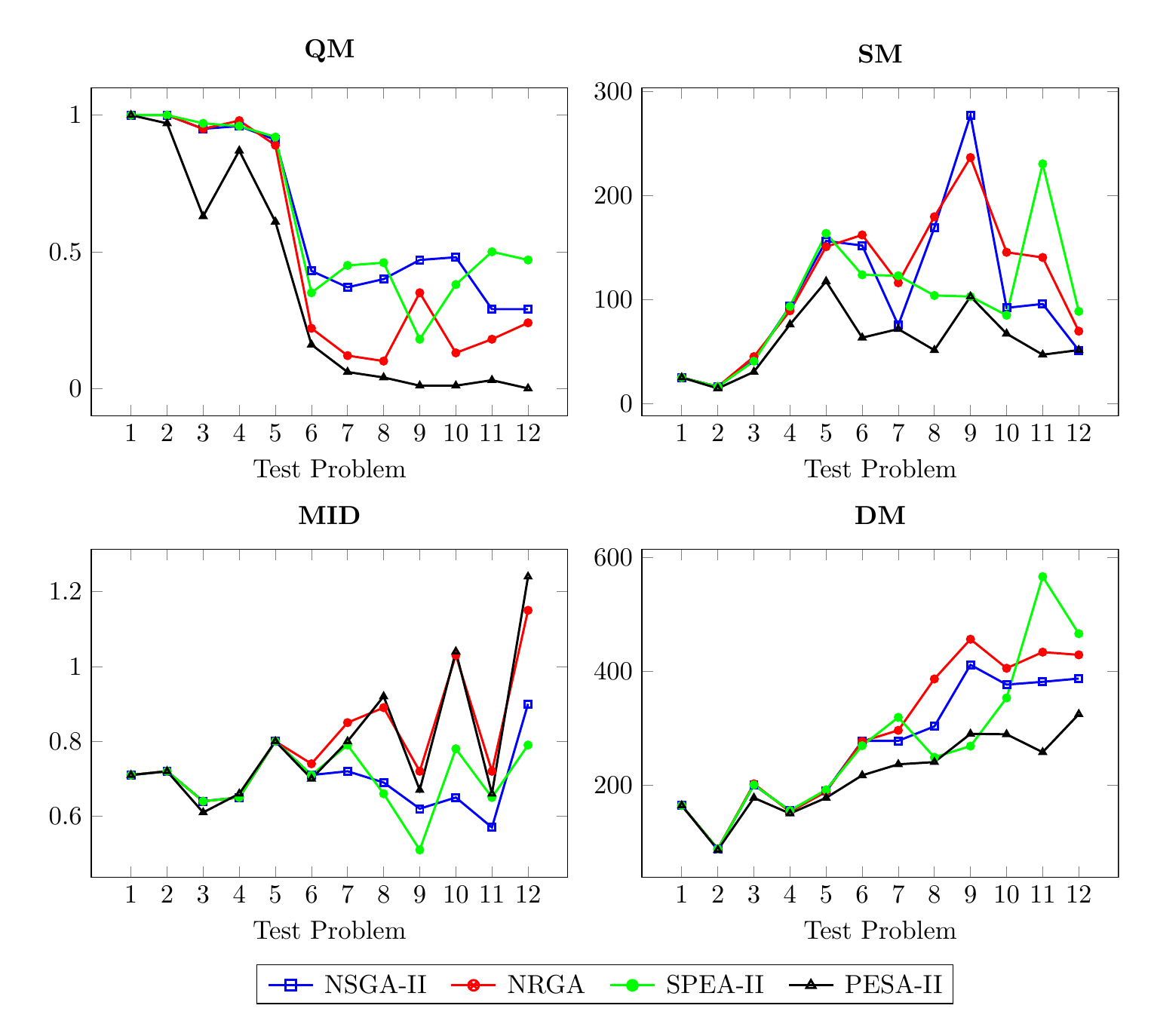}
\caption{Average results obtained by different evolutionary algorithms for different metrics }
\label{evol-compare}
\end{figure}

\begin{figure*}
\begin{multicols}{2}
    \includegraphics[width=0.95\linewidth]{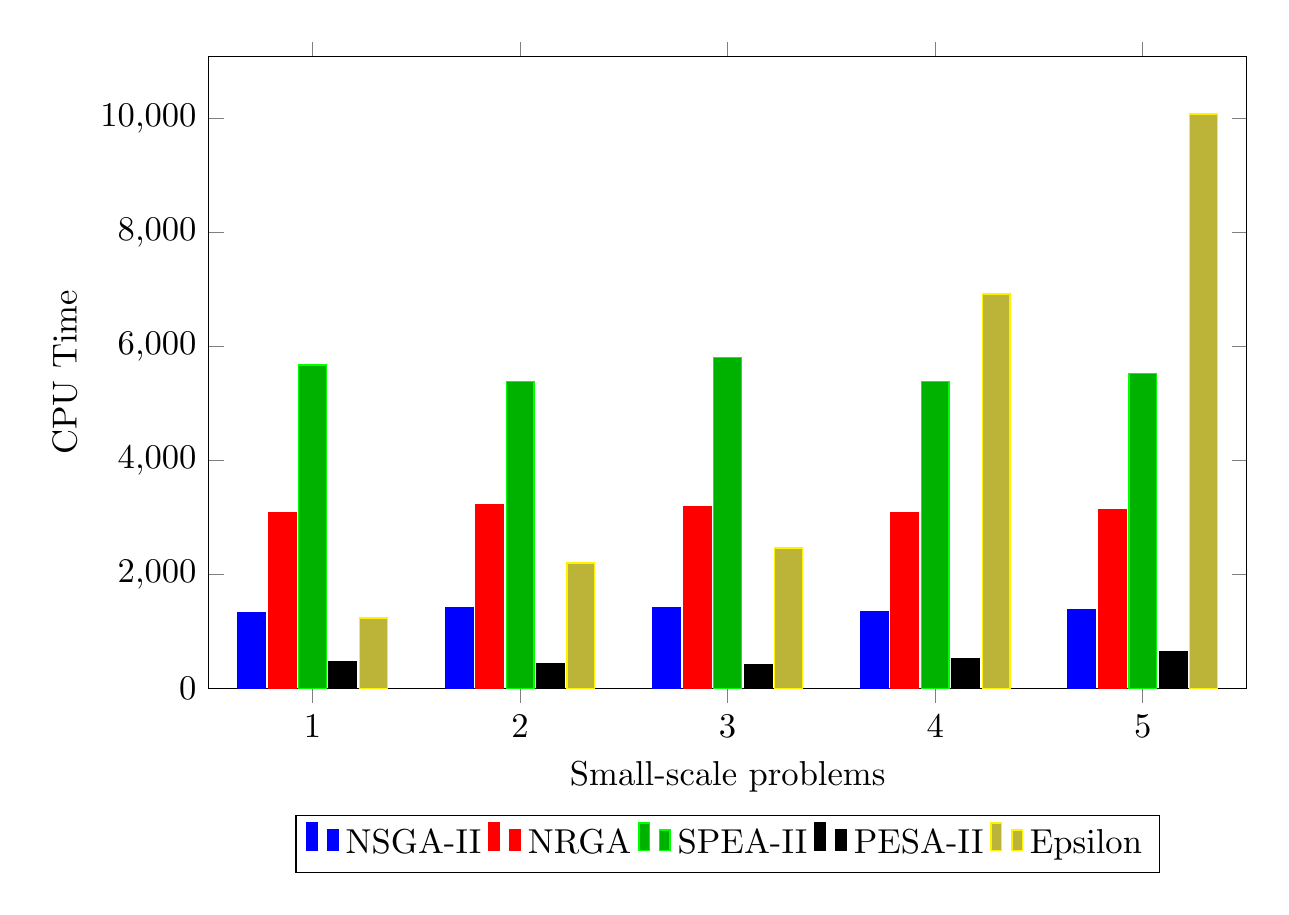}\par 
    \hspace{-12mm}
    \includegraphics[width=1.2\linewidth]{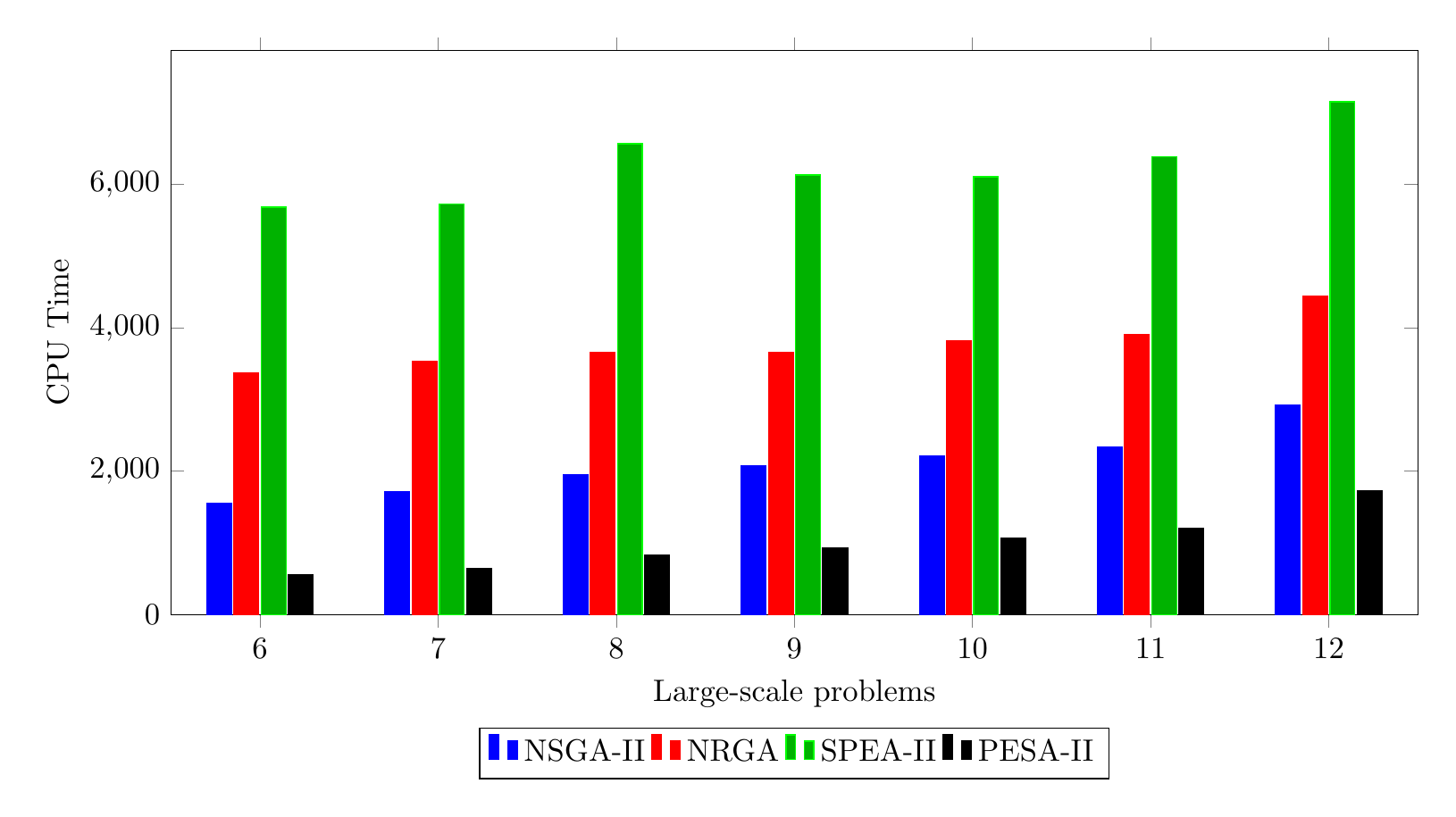}
    \end{multicols}
    \vspace*{-3mm}
    \captionsetup{justification=centering}
    \caption{Run-time comparison for evolutionary algorithms and $A\varepsilon$-c method for small-scale test problems (left), and between evolutionary algorithms for large-scale test problems (right) }
\label{time_cpu}
\end{figure*}

To evaluate the performance of the algorithms, we compare the results obtained by four algorithms with the exact methods for the small-scale test problems. Table \ref{small_scales} represents the comparison between the evolutionary algorithms considering comparison metrics with A$\varepsilon$-c for small-scale test problems for five different runs. The numerical results indicate that evolutionary algorithms could produce optimal or close to optimal solutions when compared to the exact method in small-scale problems. In addition, we illustrated the Pareto solutions obtained from NSGA-II and A$\varepsilon$-c algorithms in Figure \ref{pareto_results}. The figure shows that except for one solution, the NSGA-II was able to find all the solutions found by the exact method for the test problem 5. Also, the NSGA-II could find an additional solution where the A$\varepsilon$-c was not able to obtain due to reaching the stipulated time-limit (3 hours).
\begin{table} 
	\caption{\hspace{2.2cm}Generated test problems}
	\centering
	\begin{tabular}{|c| c c c c |} 
		\hline
		\multirow{1}{*}{Test No.}& $|k|$ & $|I|$ & $|V^{in}|$&$|V^{out}|$ \\
		\hline
		\multirow{1}{*}{1} & 2 & 4 & 3 &3 \\
  	    \multirow{1}{*}{2} & 2 & 4 & 4 &3  \\
  	    \multirow{1}{*}{3} & 2 & 4 & 3 &4 \\
  	    \multirow{1}{*}{4} &3	&	5	&	3	&	3	\\
  	    \multirow{1}{*}{5} &3	&	5	&	4	&	4	 \\
  	    \hline
  	    \multirow{1}{*}{6} &3	&	7	&	3	&	3	 \\
  	    \multirow{1}{*}{7} &4	&	10	&	5	&	5	\\
  	    \multirow{1}{*}{8} &5	&	15	&	7	&	7	\\
  	    \multirow{1}{*}{9} &6	&	20	&	9	&	9	\\
  	    \multirow{1}{*}{10}&7	&	25	&	11	&	11	\\
  	    \multirow{1}{*}{11}&8	&	30	&	13	&	13	\\
  	    \multirow{1}{*}{12}&10	&	50	&	15	&	15\\
		\hline
	\end{tabular}
	\label{instances}
\end{table}
Table \ref{large-scale-problems} presents the numerical results for the large-scale test problems and Figure \ref{evol-compare} compares the average results obtained by different evolutionary algorithms considering multiple metrics for all test problems. It should be noted that the termination criterion for all algorithms is reaching to $3\times 10^5$ for the number of function evaluations. In addition, we compare the solution methods in terms of running time provided in Figure \ref{time_cpu}. The results show that in all instances, the PESA-II algorithm significantly outperformed other algorithm by. Other than PESA-II, the NSGA-II performed faster than other algorithms in all test problems. Also, The A$\varepsilon$-c could only outperformed NRGA and SPEA in the first three test problems. Table \ref{large-scale-problems} presents the numerical results for the large-scale test problems. Figure \ref{evol-compare} compares the quality of the solutions obtained by different evolutionary algorithms considering multiple metrics for all test problems.

\begin{table*}
	\caption{\hspace{2.2cm} Numerical results for small-scale problems}
	\centering
	\begin{adjustbox}{width=1.1\textwidth,center=\textwidth}
	\begin{tabular}{|c|c c c c c| c c c c c| c c c c c| c c c c c|} 
		\hline
		\multirow{3}{*}{}&
		\multicolumn{20}{c|}{Metrics \& Algorithms}\\
	     \cline{2-21}
		 \multirow{1}{*}{Tests problems No.}& \multicolumn{5}{c|}{QM}& \multicolumn{5}{c|}{SM}& \multicolumn{5}{c|}{MID}& \multicolumn{5}{c|}{DM}\\
		\cline{2-6} \cline{7-11} \cline{12-16} \cline{17-21}
		 & NSGA-II& NRGA & SPEA-II & PESA-II &  A$\varepsilon$c & NSGA-II& NRGA & SPEA-II & SPEA-II &A$\varepsilon$c & NSGA-II& NRGA & SPEA-II & PESA-II & A$\varepsilon$c & NSGA-II& NRGA & SPEA-II  & PESA-II &A$\varepsilon$-c \\
		 \hline
		 \multirow{5}{*}{1} & 1.00&1.00&1.00&1.00&1.00&24.83&24.83&24.83&24.83&24.83&0.71&0.71&0.71&0.71&0.71& 164.90&164.90&164.90&164.90&164.90\\
		 &1	&	1	&	1	&	1	&	1	&	24.83	&	24.83	&	24.83	&	24.83	&	24.83	&	0.71	&	0.71	&	0.71	&	0.71	&	0.71	&	164.9	&	164.9	&	164.9	&	164.9	&	164.9\\
		 &1	&	1	&	1	&	1	&	1	&	24.83	&	24.83	&	24.83	&	24.83	&	24.83	&	0.71	&	0.71	&	0.71	&	0.71	&	0.71	&	164.9	&	164.9	&	164.9	&	164.9	&	164.9\\
		 &1	&	1	&	1	&	1	&	1	&	24.83	&	24.83	&	24.83	&	24.83	&	24.83	&	0.71	&	0.71	&	0.71	&	0.71	&	0.71	&	164.9	&	164.9	&	164.9	&	164.9	&	164.9\\
		 &1	&	1	&	1	&	1	&	1	&	24.83	&	24.83	&	24.83	&	24.83	&	24.83	&	0.71	&	0.71	&	0.71	&	0.71	&	0.71	&	164.9	&	164.9	&	164.9	&	164.9	&	164.9\\
                      
		\hline
		\multirow{5}{*}{2} &   1	&	1	&	1	&	1	&	1	&	16.07	&	16.07	&	16.07	&	16.07	&	16.07	&	0.72	&	0.72	&	0.72	&	0.72	&	0.72	&	88.63	&	88.63	&	88.63	&	88.63	&	88.63\\
				&	1	&	1	&	1	&	1	&	1	&	16.07	&	16.07	&	16.07	&	16.07	&	16.07	&	0.72	&	0.72	&	0.72	&	0.72	&	0.72	&	88.63	&	88.63	&	88.63	&	88.63	&	88.63\\
				&	1	&	1	&	1	&	1	&	1	&	16.07	&	16.07	&	16.07	&	16.07	&	16.07	&	0.72	&	0.72	&	0.72	&	0.72	&	0.72	&	88.63	&	88.63	&	88.63	&	88.63	&	88.63\\
				&	1	&	1	&	1	&	1	&	1	&	16.07	&	16.07	&	16.07	&	16.07	&	16.07	&	0.72	&	0.72	&	0.72	&	0.72	&	0.72	&	88.63	&	88.63	&	88.63	&	88.63	&	88.63\\
				&	1	&	1	&	1	&	0.86	&	1	&	16.07	&	16.07	&	16.07	&	7.51	&	16.07	&	0.72	&	0.72	&	0.72	&	0.68	&	0.72	&	88.63	&	88.63	&	88.63	&	80.08	&	88.63\\

		\hline
		\multirow{5}{*}{3} & 1	&	1	&	0.88	&	0.5	&	1	&	43.07	&	43.07	&	30.36	&	15.92	&	43.07	&	0.63	&	0.63	&	0.67	&	0.63	&	0.63	&	203.1	&	203.1	&	191.8	&	171.4	&	203.1\\
				&	0.88	&	0.88	&	1	&	0.75	&	1	&	47.7	&	47.7	&	43.07	&	47.82	&	43.07	&	0.64	&	0.64	&	0.63	&	0.64	&	0.63	&	202.4	&	202.4	&	203.1	&	203.9	&	203.1\\
				&	1	&	1	&	1	&	0.25	&	1	&	43.07	&	43.07	&	43.07	&	13.25	&	43.07	&	0.63	&	0.63	&	0.63	&	0.45	&	0.63	&	203.1	&	203.1	&	203.1	&	122	&	203.1\\
				&	1	&	0.88	&	1	&	0.75	&	1	&	43.07	&	47.7	&	43.07	&	28.94	&	43.07	&	0.63	&	0.64	&	0.63	&	0.67	&	0.63	&	203.1	&	202.4	&	203.1	&	193.5	&	203.1\\
				&	0.88	&	1	&	1	&	0.88	&	1	&	30.36	&	43.07	&	43.07	&	45.49	&	43.07	&	0.67	&	0.63	&	0.63	&	0.64	&	0.63	&	191.8	&	203.1	&	203.1	&	200.5	&	203.1\\

		\hline
		\multirow{5}{*}{4} &  1	&	1	&	0.91	&	0.73	&	0.91	&	92.47	&	92.47	&	93.52	&	50.7	&	91.51	&	0.65	&	0.65	&	0.65	&	0.67	&	0.69	&	156.1	&	156.1	&	155.8	&	142.7	&	158.2\\
				&	0.91	&	0.91	&	0.91	&	0.82	&	0.91	&	93.52	&	75.57	&	95.74	&	50.02	&	91.51	&	0.65	&	0.62	&	0.65	&	0.67	&	0.69	&	155.8	&	145.3	&	155.2	&	143	&	158.2\\
				&	1	&	1	&	1	&	1	&	0.91	&	92.47	&	92.47	&	92.47	&	92.47	&	91.51	&	0.65	&	0.65	&	0.65	&	0.65	&	0.69	&	156.1	&	156.1	&	156.1	&	156.1	&	158.2\\
				&	1	&	1	&	1	&	0.91	&	0.91	&	92.47	&	92.47	&	92.47	&	93.52	&	91.51	&	0.65	&	0.65	&	0.65	&	0.65	&	0.69	&	156.1	&	156.1	&	156.1	&	155.8	&	158.2\\
				&	0.91	&	1	&	1	&	0.91	&	0.91	&	96.71	&	92.47	&	92.47	&	92.57	&	91.51	&	0.65	&	0.65	&	0.65	&	0.65	&	0.69	&	155	&	156.1	&	156.1	&	156	&	158.2\\

		\hline
		\multirow{5}{*}{5} & 0.93	&	0.93	&	1	&	0.67	&	0.87	&	153.6	&	161.9	&	188.4	&	119.9	&	153.5	&	0.8	&	0.8	&	0.79	&	0.79	&	0.78	&	191.5	&	191.7	&	197.3	&	178.7	&	190.2\\
				&	0.93	&	0.93	&	1	&	0.67	&	0.87	&	161.9	&	153.6	&	188.4	&	91.75	&	153.5	&	0.8	&	0.8	&	0.79	&	0.8	&	0.78	&	191.7	&	191.5	&	197.3	&	172.6	&	190.2\\
				&	0.93	&	0.87	&	0.73	&	0.73	&	0.87	&	153.6	&	184.6	&	125.4	&	96.55	&	153.5	&	0.8	&	0.79	&	0.82	&	0.78	&	0.78	&	191.5	&	192.7	&	183.8	&	170.5	&	190.2\\
				&	0.8	&	0.8	&	0.93	&	0.47	&	0.87	&	150.4	&	99.99	&	153.6	&	97.49	&	153.5	&	0.8	&	0.79	&	0.8	&	0.81	&	0.78	&	186.7	&	176.9	&	191.5	&	176.6	&	190.2\\
				&	0.93	&	0.93	&	0.93	&	0.53	&	0.87	&	161.9	&	153.6	&	161.9	&	181.2	&	153.5	&	0.8	&	0.8	&	0.8	&	0.81	&	0.78	&	191.7	&	191.5	&	191.7	&	192.6	&	190.2\\

		\hline
	\end{tabular}
	\end{adjustbox}
	\label{small_scales}
\end{table*}
%%%%%%%%%%%%%%%%%%%%%%%%%%%%%%%%%%%%%%%%%%%%%%%%%%%%%%%%%%%%
\begin{table*}
	\caption{\hspace{2.2cm} Numerical results for large-scale problems}
	\centering
	\begin{adjustbox}{width=1.1\textwidth,center=\textwidth}
	\begin{tabular}{|c|c c c c c| c c c c c| c c c c c| c c c c c|} 
		\hline
		\multirow{3}{*}{}&
		\multicolumn{20}{c|}{Metrics \& Algorithms}\\
	     \cline{2-21}
		 \multirow{1}{*}{Tests problems}& \multicolumn{5}{c|}{QM}& \multicolumn{5}{c|}{SM}& \multicolumn{5}{c|}{MID}& \multicolumn{5}{c|}{DM}\\
		\cline{2-6} \cline{7-11} \cline{12-16} \cline{17-21}
		 & NSGA-II& NRGA & SPEA-II & PESA-II &  $A\varepsilon$-c  & NSGA-II& NRGA & SPEA-II & SPEA-II &$A\varepsilon$-c  & NSGA-II& NRGA & SPEA-II & PESA-II &$A\varepsilon$-c  & NSGA-II& NRGA & SPEA-II  & PESA-II &$A\varepsilon$-c \\
		 		\multirow{5}{*}{6} &0.43	&	0.39	&	0.17	&	0	&	-	&	90.26	&	169.5	&	115.3	&	40.4	&	-	&	0.68	&	0.71	&	0.72	&	0.67	&	-	&	243.9	&	295.7	&	281.9	&	182.2	&	-\\
				&	0.6	&	0.05	&	0.3	&	0.1	&	-	&	187.8	&	181.5	&	160.7	&	86.96	&	-	&	0.69	&	0.77	&	0.69	&	0.8	&	-	&	288.5	&	288.1	&	290.6	&	260.7	&	-\\
				&	0.18	&	0.41	&	0.47	&	0.24	&	-	&	105.8	&	87.45	&	111.3	&	32.99	&	-	&	0.71	&	0.69	&	0.74	&	0.69	&	-	&	285.5	&	241.3	&	260.2	&	207.6	&	-\\
				&	0.48	&	0.24	&	0.52	&	0	&	-	&	116.6	&	267.3	&	133.1	&	31.46	&	-	&	0.74	&	0.7	&	0.7	&	0.65	&	-	&	266.6	&	314.3	&	268.3	&	161.3	&	-\\
				&	0.47	&	0	&	0.26	&	0.47	&	-	&	258.9	&	104.7	&	98.24	&	124.1	&	-	&	0.72	&	0.84	&	0.7	&	0.71	&	-	&	305.1	&	247.4	&	247.4	&	276.6	&	-\\

                    \hline
		 
		 \hline
		 \multirow{5}{*}{7} &0.36	&	0.14	&	0.5	&	0	&	-	&	43.28	&	92.83	&	104	&	59.89	&	-	&	0.79	&	1.19	&	0.86	&	1	&	-	&	255.2	&	269	&	280.6	&	200.6	&	-\\
				&	0.52	&	0.05	&	0.43	&	0	&	-	&	103.3	&	171	&	266.9	&	89.8	&	-	&	0.77	&	0.78	&	0.77	&	0.71	&	-	&	321.6	&	342.6	&	413.6	&	233.9	&	-\\
				&	0.23	&	0.08	&	0.54	&	0.15	&	-	&	78.68	&	133.4	&	70.23	&	25.39	&	-	&	0.74	&	0.79	&	0.81	&	0.73	&	-	&	288.2	&	324	&	326.1	&	162.3	&	-\\
				&	0.31	&	0.15	&	0.54	&	0	&	-	&	32	&	116.8	&	113.9	&	87.48	&	-	&	0.64	&	0.74	&	0.75	&	0.8	&	-	&	235.5	&	275.2	&	309	&	276.3	&	-\\
				&	0.43	&	0.19	&	0.24	&	0.14	&	-	&	118.7	&	66.15	&	59.13	&	94.96	&	-	&	0.66	&	0.74	&	0.77	&	0.75	&	-	&	290.3	&	272.1	&	267.8	&	312.2	&	-\\

 		\hline
		 \multirow{5}{*}{8} &1	&	0	&	0	&	0	&	-	&	35.71	&	95.57	&	78.17	&	109.1	&	-	&	0.76	&	1.16	&	1.05	&	1.1	&	-	&	192.9	&	335.3	&	255.7	&	230.4	&	-\\
				&	0	&	0	&	1	&	0	&	-	&	57.76	&	157.2	&	128.3	&	6.66	&	-	&	1.04	&	1	&	0.73	&	1.68	&	-	&	261.5	&	365	&	250.5	&	184.6	&	-\\
				&	0.26	&	0.42	&	0.32	&	0	&	-	&	307.7	&	226.8	&	22.32	&	50.08	&	-	&	0.39	&	0.66	&	0.41	&	0.53	&	-	&	339.3	&	444.1	&	178.7	&	224.8	&	-\\
				&	0.5	&	0.06	&	0.44	&	0	&	-	&	232.1	&	51.04	&	163.7	&	41.42	&	-	&	0.46	&	0.65	&	0.42	&	0.61	&	-	&	369	&	316.4	&	290.4	&	291.4	&	-\\
				&	0.25	&	0	&	0.55	&	0.2	&	-	&	212.7	&	367.5	&	126.7	&	48.54	&	-	&	0.77	&	0.97	&	0.67	&	0.69	&	-	&	354.8	&	472.7	&	271.1	&	273.1	&	-\\

		\hline
		 \multirow{5}{*}{9} &0.38	&	0.58	&	0.04	&	0	&	-	&	251.9	&	256.9	&	47.62	&	496.7	&	-	&	0.45	&	0.56	&	0.4	&	0.66	&	-	&	438.3	&	498	&	157.1	&	673.3	&	-\\
				&	0.21	&	0.57	&	0.21	&	0	&	-	&	100.5	&	121.2	&	72.69	&	32.29	&	-	&	0.79	&	0.65	&	0.52	&	0.79	&	-	&	365	&	340	&	251.4	&	251.8	&	-\\
				&	0.61	&	0.18	&	0.18	&	0.04	&	-	&	361.7	&	299	&	248.4	&	20.3	&	-	&	0.51	&	0.67	&	0.53	&	0.58	&	-	&	388.2	&	486.5	&	399.8	&	256.2	&	-\\
				&	0.65	&	0	&	0.35	&	0	&	-	&	210.2	&	301.6	&	47.94	&	52.99	&	-	&	0.68	&	1.13	&	0.5	&	0.73	&	-	&	352.4	&	545.3	&	224.4	&	196.1	&	-\\
				&	0.48	&	0.4	&	0.12	&	0	&	-	&	462.8	&	204.3	&	97.77	&	10.31	&	-	&	0.65	&	0.58	&	0.61	&	0.61	&	-	&	513.4	&	412.1	&	312.2	&	74.76	&	-\\

		\hline
		 \multirow{5}{*}{10} &0.25	&	0.15	&	0.6	&	0	&	-	&	69.47	&	146.1	&	84.63	&	76.16	&	-	&	0.78	&	0.92	&	0.65	&	1.18	&	-	&	366.2	&	429.8	&	294.7	&	288.6	&	-\\
				&	0.33	&	0.33	&	0.29	&	0.05	&	-	&	136	&	262.4	&	49.55	&	5.45	&	-	&	0.61	&	0.64	&	0.6	&	0.62	&	-	&	391.5	&	527.4	&	355	&	184.1	&	-\\
				&	0.61	&	0	&	0.39	&	0	&	-	&	129.8	&	230.7	&	111.7	&	162.2	&	-	&	0.73	&	1.56	&	1.06	&	1.82	&	-	&	382.3	&	510.6	&	380.4	&	402	&	-\\
				&	0.36	&	0	&	0.64	&	0	&	-	&	50.73	&	20.45	&	81.2	&	89.42	&	-	&	0.78	&	1.36	&	0.78	&	1.29	&	-	&	442.6	&	224.6	&	325.3	&	517.2	&	-\\
				&	0.83	&	0.17	&	0	&	0	&	-	&	73.36	&	67.64	&	96.82	&	1.88	&	-	&	0.39	&	0.65	&	0.81	&	0.3	&	-	&	301.7	&	335.3	&	411.6	&	56.43	&	-\\

		\hline
		 \multirow{5}{*}{11} &0.29	&	0.21	&	0.46	&	0.04	&	-	&	20.98	&	220.6	&	441.1	&	45.09	&	-	&	0.42	&	0.69	&	0.6	&	0.54	&	-	&	238.6	&	439.7	&	721	&	264.8	&	-\\
				&	0.26	&	0.37	&	0.37	&	0	&	-	&	99.84	&	62.56	&	212.1	&	54.85	&	-	&	0.51	&	0.61	&	0.71	&	0.68	&	-	&	387.8	&	402.8	&	603.4	&	274.9	&	-\\
				&	0.14	&	0.14	&	0.68	&	0.05	&	-	&	59.48	&	27.3	&	193.9	&	89.62	&	-	&	0.61	&	0.7	&	0.59	&	0.63	&	-	&	408.6	&	249.8	&	501.9	&	304.7	&	-\\
				&	0.47	&	0.18	&	0.35	&	0	&	-	&	75.53	&	116.4	&	231.6	&	15.15	&	-	&	0.77	&	1.05	&	0.85	&	1.06	&	-	&	378	&	566.1	&	550.3	&	205.5	&	-\\
				&	0.29	&	0	&	0.64	&	0.07	&	-	&	222.5	&	275.2	&	73.75	&	29.32	&	-	&	0.52	&	0.56	&	0.52	&	0.38	&	-	&	495.6	&	510.5	&	454.3	&	241.1	&	-\\

		\hline
		 \multirow{5}{*}{12} &0	&	0.2	&	0.8	&	0	&	-	&	39.72	&	100.8	&	28.73	&	174.3	&	-	&	1.22	&	1.17	&	0.69	&	1.83	&	-	&	327.9	&	477.9	&	254.8	&	636.3	&	-\\
				&	0.53	&	0	&	0.47	&	0	&	-	&	24.89	&	56.48	&	64.42	&	38.78	&	-	&	0.62	&	0.71	&	0.75	&	0.86	&	-	&	379.3	&	404.1	&	494.7	&	309.4	&	-\\
				&	0	&	0.75	&	0.25	&	0	&	-	&	59.24	&	49.93	&	192.9	&	8.4	&	-	&	0.79	&	0.96	&	0.92	&	1.14	&	-	&	417.2	&	406.9	&	623.4	&	130.3	&	-\\
				&	0.7	&	0	&	0.3	&	0	&	-	&	55.05	&	66.35	&	59.61	&	20.74	&	-	&	1.14	&	1.75	&	0.89	&	1.38	&	-	&	413.1	&	456.8	&	443.3	&	309.3	&	-\\
				&	0.2	&	0.27	&	0.53	&	0	&	-	&	75.61	&	73.62	&	96.88	&	13.99	&	-	&	0.76	&	1.14	&	0.72	&	0.98	&	-	&	399.5	&	399.2	&	414.5	&	239.2	&	-\\

		\hline

	\end{tabular}
	\label{large-scale-problems}
\end{adjustbox}
	\label{Decomposition_2}
\end{table*}

% \section{Results}

% This section may be divided by subheadings. It should provide a concise and precise description of the experimental results, their interpretation as well as the experimental conclusions that can be drawn.
% \begin{quote}
% This section may be divided by subheadings. It should provide a concise and precise description of the experimental results, their interpretation as well as the experimental conclusions that can be drawn.
% \end{quote}

%%%%%%%%%%%%%%%%%%%%%%%%%%%%%%%%%%%%%%%%%%
\subsubsection{Statistical Analysis}
In this section, we investigate whether there is any statistically significant difference across the proposed methods used for each defined metrics in the methodology section. We use the non-parametric Kruskal-Wallis test \cite{kruskal1952use} since the sample size (here test problems) is small, the normality assumption is not met, and the number of methods is larger than two. We use 0.05 as a threshold for the p-value of the test to decide whether to reject the null hypothesis, indicating that there are at least two methods whose distribution of metrics are different from each other for large scale problems. We use SPSS software version 26 to first examine the QM distribution. Applying Kruskal-Wallis test, the resulting p-value is $<0.0005$ and the null hypothesis is rejected. 

\begin{figure*}
\centering
\captionsetup{justification=centering}
\includegraphics[scale=1]{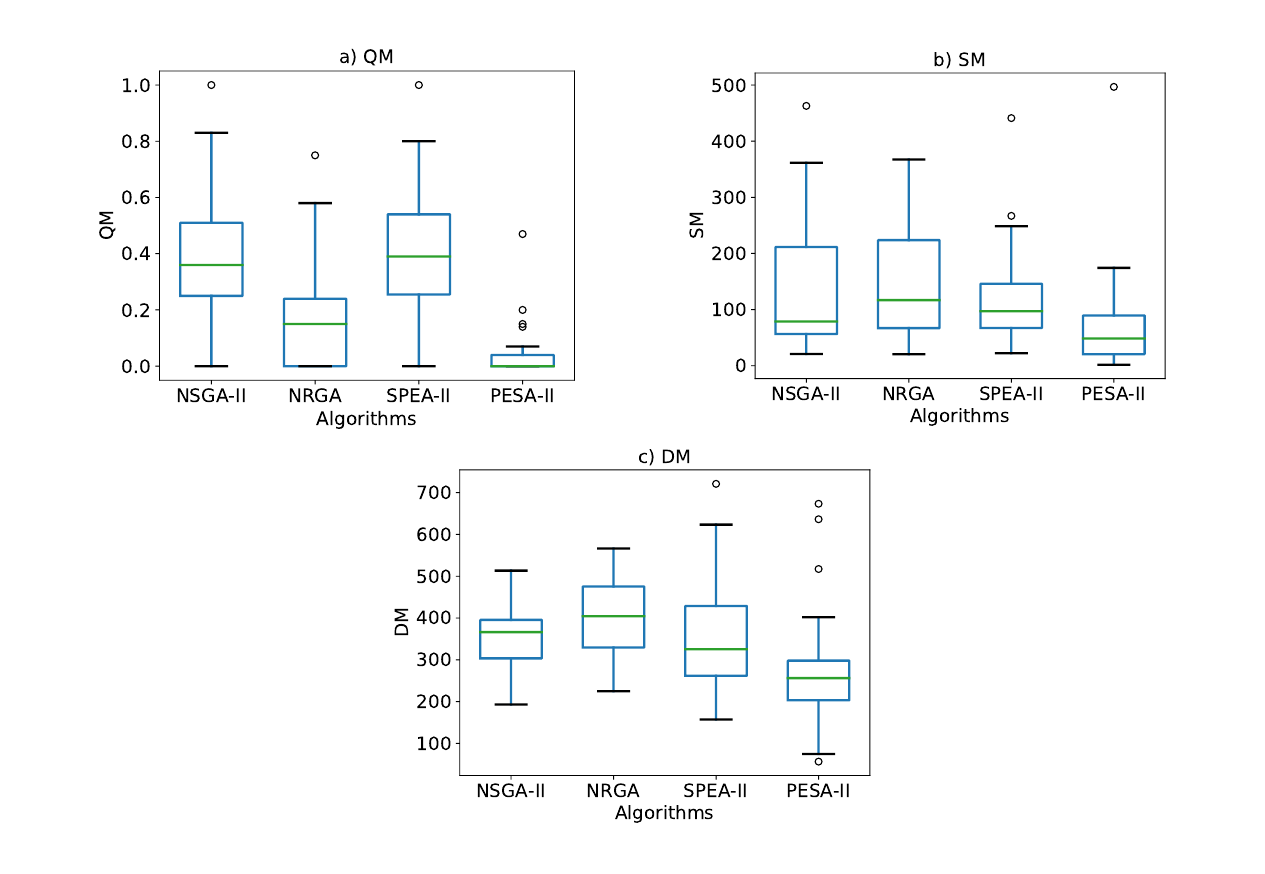}\hfill
\caption{Variability in each comparison metric across different algorithms}
\label{fig:boxplot}
\end{figure*}

To find the exact methods whose QM metric is different from other algorithms, we run a pairwise comparison test and show the results in Table \ref{pairwise-QM}. Each row tests the null hypothesis that Sample 1 and Sample 2 distributions are the same. Asymptotic significances (2-sided tests) are displayed in the fifth column. Significance values have been adjusted by the Bonferroni correction for multiple testing, shown in the last column.
Methods NSGA-II and PESA-II with p-value of 0.001 and methods SPEA-II and PESA-II with the same p-value are found to be statistically significantly different from each other in terms of QM values (Figure \ref{fig:boxplot}-a). 

For SM, the p-value of Kruskal-Wallis test is 0.016, concluding that the null hypothesis is rejected. The pairwise comparison tests in Table \ref{pairwise-SM} shows that the distribution of SM for methods PESA-II and NRGA are statistically significantly different from each other (p-value of 0.012) (Figure \ref{fig:boxplot}-b). Since a lower value for SM is more preferable, we conclude the PESA-II method performs better than other methods regarding this specific metric.

While applying Kruskal-Wallis test on MID results in a p-value of 0.064 and indicates that we cannot reject the null hypothesis, the same p-value for DM is 0.031. Using pairwise comparisons tests (Table \ref{pairwise-DM}), we find that methods PESA-II and NRGA are statistically significantly different from each other in terms of DM values with p-value of 0.023. Since a higher value for DM is more preferable, methods NRGA and NSGA-II are performing better than other methods in terms of DM metric (Figure \ref{fig:boxplot}-c).
\begin{table*}
	\caption{\hspace{2.2cm} Pairwise comparisons of methods for QM}
	\centering
	\begin{tabular}{|c| c c c c c|} 
		\hline
		\multirow{1}{*}{Sample 1 - Sample 2}& Test Statistics & Std. Error & Std. Test Statistic & Sig.& Adj.Sig.\\
		\hline
		\multirow{1}{*}{PESA-II - NRGA} & 7.000 & 4.394 & 1.593 &0.111&0.667  \\
  	    \multirow{1}{*}{PESA-II - NSGA-II} & 16.357& 4.394 & 3.723 &0.000&\textcolor{red}{0.001} \\
  	    \multirow{1}{*}{PESA-II - SPEA-II} &-16.929	&	4.394&	-3.853	&	0.000& \textcolor{red}{0.001}\\
  	    \multirow{1}{*}{NRGA - NSGA-II} & -9.357	&4.394	&-2.130&	0.033&0.199\\
  	    \multirow{1}{*}{NRGA - SPEA-II} &-9.929	&	4.394	&	-2.260	&	0.024&0.143\\
  	    \multirow{1}{*}{NSGA-II - SPEA-II} & -0.571	&	4.394	&	-0.130	&0.897&1.000\\
  	    \hline
	\end{tabular}
	\label{pairwise-QM}
\end{table*}
\begin{table*}
	\caption{\hspace{2.2cm} Pairwise comparisons of methods for SM}
	\centering
	\begin{tabular}{|c| c c c c c|} 
		\hline
		\multirow{1}{*}{Sample 1 - Sample 2}& Test Statistics & Std. Error & Std. Test Statistic & Sig.& Adj.Sig.\\
		\hline
		\multirow{1}{*}{PESA-II - NRGA} & 13.571 & 4.397 &3.087 &0.002&\textcolor{red}{0.012} \\
  	    \multirow{1}{*}{PESA-II - NSGA-II} & 9.429& 4.397 & 2.144 &0.032&0.192\\
  	    \multirow{1}{*}{PESA-II - SPEA-II} &-9.857	&	4.397&	-2.242	&	0.025& 0.150\\
  	    \multirow{1}{*}{NRGA - NSGA-II} & 4.143	&4.397	&0.942&	0.346&1.000\\
  	    \multirow{1}{*}{NRGA - SPEA-II} &3.714	&	4.397	&0.845	&0.398&1.000\\
  	    \multirow{1}{*}{NSGA-II - SPEA-II} & -0.429	&	4.397	&	-0.097	&0.922&1.000\\
  	    \hline
	\end{tabular}
	\label{pairwise-SM}
\end{table*}
\begin{table*}
	\caption{\hspace{2.2cm} Pairwise comparisons of methods for DM}
	\centering
	\begin{tabular}{|c| c c c c c|} 
		\hline
		\multirow{1}{*}{Sample 1 - Sample 2}& Test Statistics & Std. Error & Std. Test Statistic & Sig.& Adj.Sig.\\
		\hline
		\multirow{1}{*}{PESA-II - NRGA} & 12.714& 4.397 & 2.892 &0.004&\textcolor{red}{0.023} \\
  	    \multirow{1}{*}{PESA-II - NSGA-II} & 9.143& 4.397 & 2.079 &.038&0.226\\
  	    \multirow{1}{*}{PESA-II - SPEA-II} &-7.571	&	4.397&	-1.722	&	0.085& 0.510\\
  	    \multirow{1}{*}{NRGA - NSGA-II} & 3.571	&4.397	&0.812&	0.417&1.000\\
  	    \multirow{1}{*}{NRGA - SPEA-II} &5.142&	4.397	&1.170	&	0.242&1.000\\
  	    \multirow{1}{*}{NSGA-II - SPEA-II} & 1.571	&	4.397	&	0.357	&0.721&1.000\\
  	    \hline
	\end{tabular}
	\label{pairwise-DM}
\end{table*}

% \begin{figure}
%     %  \centering
%      \begin{subfigure}
%          \centering
%          \includegraphics[width=0.4\textwidth]{DM.pdf}
%          \caption{(a)}
%          \label{dm}
%      \end{subfigure}
%      \hfill
%      \begin{subfigure}
%         %  \centering
%          \includegraphics[width=0.4\textwidth]{QM.pdf}
%          \caption{(b)}
%          \label{qm}
%      \end{subfigure}
%      \hfill
%      \begin{subfigure}
%         %  \centering
%          \includegraphics[width=0.4\textwidth]{SM.pdf}
%          \caption{(c)}
%          \label{sm}
%      \end{subfigure}
%         \caption{Three simple graphs}
%         \label{fig:boxplot}
% \end{figure}
\subsection{Sensitivity Analysis}
\begin{figure}
\centering
\includegraphics[scale=0.7]{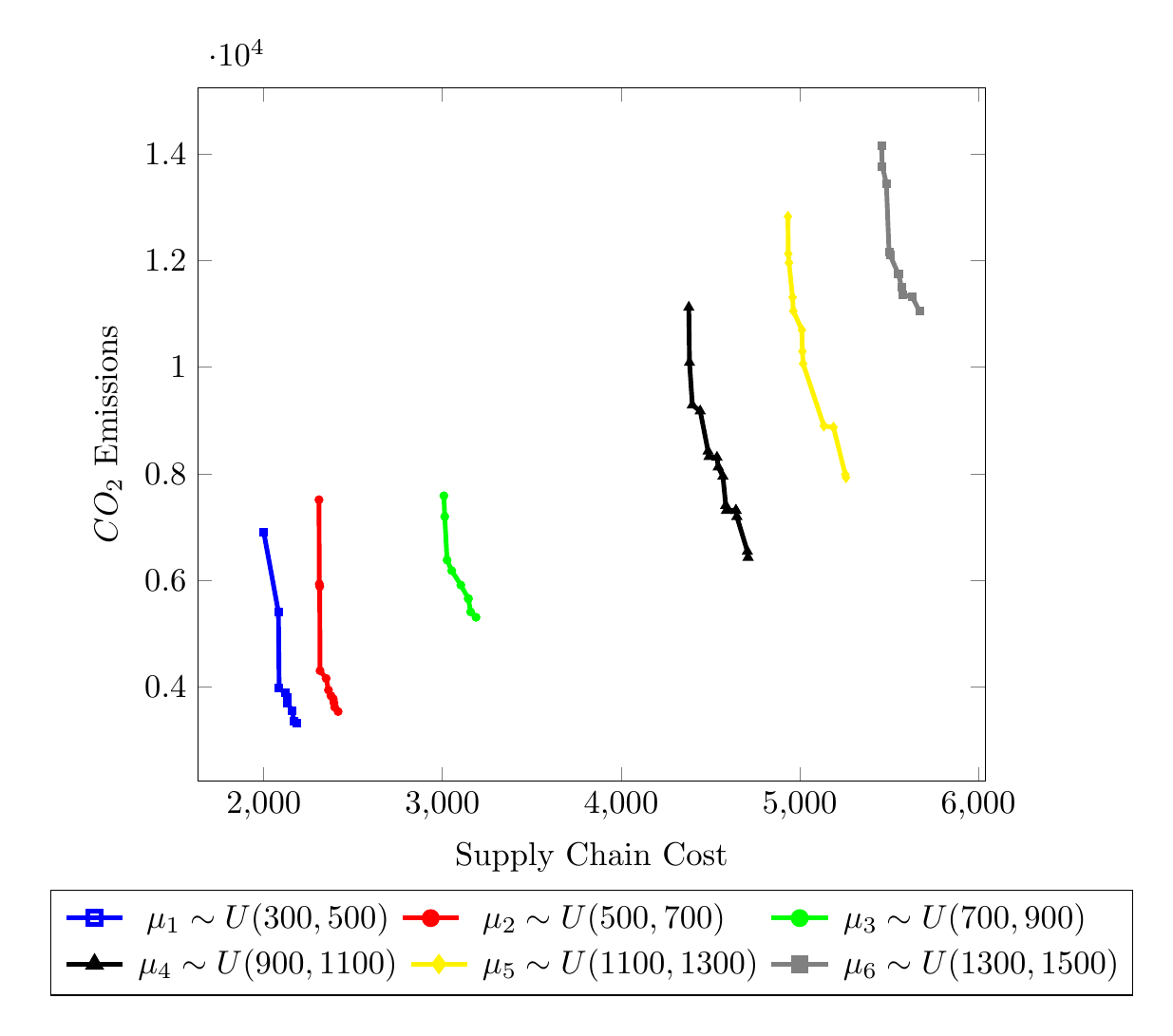}
\captionsetup{justification=centering}
\caption{The Pareto solutions obtained under different scenarios for the mean of annual demand.}
\label{sens_1}
\end{figure}

% \subsubsection{Effect of Demand Parameter}
We analyzed the effect of changes in the mean of annual demand.  We considered six scenarios where we randomly draw six numbers from the scenario intervals. As it is shown in Figure \ref{sens_1}, the proposed model is quite sensitive to the demand parameter. At each scenario, multiple solutions are provided by solving the model. These sets of solutions enable the decision makers to significantly decrease the amount of $CO_{2}$ emission with a small sacrifice in the supply chain cost. This can help companies to move towards the sustainable supply chain which will help them to mitigate the adverse effects of green house gas emission and improve their competitiveness in the market. 
In addition, it can be observed that with an increase in the demand parameter we see an increase in the supply chain cost as well as the $CO_{2}$ emissions. This is mainly due to the increase in transportation activities and the need for establishment of more DCs because of higher demand volume.
% \subsubsection{Effect of Demand on Inventory}
% We analyzed the effect demand on inventory by considering the same scenarios as the previous experiment. To better interpret these experiments, we define $\zeta_1$ and $\zeta_2$ as the weights of first and second objectives. For each scenario we considers three different sub-scenarios for objective weights as ($\zeta_1=1$, $\zeta_2=0$), ($\zeta_1=0$, $\zeta_2=1$), ($\zeta_1=0.5$, $\zeta_2=0.5$).

%====================================================================
%%%%%%%%%%%%%%%%%%%%%%%%%% Appendix %%%%%%%%%%%%%%%%%%%%%%%%%%%%%%%%%
%====================================================================
% \clearpage
\section{Conclusion}\label{con}
In this research, we proposed a new mathematical model to design a green supply chain under uncertain data parameters. The proposed model incorporates different inventory component into the location-routing decisions. Besides, we considered the carbon emissions caused by transportation activities as well as distribution centers into the proposed framework to suggest a sustainable supply chain. To efficiently solve the proposed model, we first reformulated the model to a linear mixed-integer problem and then implemented an exact method which adapts to the multi-objective problems. For large-scale instances, we proposed four evolutionary algorithms with different characteristics. Extensive numerical experiments are performed on randomly generated test problems. The performance of four evolutionary algorithms was tested by comparing it to the exact method. For the small-scale instances, evolutionary algorithms were able to obtain optimal or near-optimal solutions. In terms of run-time, PESA-II performed quicker comparing to other algorithms. We also performed a statistical analysis to compare the evolutionary algorithms in terms of the quality of the solutions. Our statistical analysis showed that methods NSGA-II and SPEA-II are performing better in terms of the QM metric. For the SM metric, the PESA method demonstrates better efficiency compared to other algorithms.
Moreover, methods NSGA-II and NRGA demonstrated better performance in regard to DM metric. To validate the model, we performed a case study, where we tested the model under different scenarios for the demand parameter. The results showed that on any demand level, the proposed model suggests a set of solutions where with a minor increase in supply chain cost, we remarkably decrease the $CO_2$ emissions. These solutions also provide various options for decision-makers to select their preferred outcome. 
Future work could include incorporating multiple scenarios into the proposed model. However, this will significantly increase the complexity of the proposed model, which requires designing efficient heuristic methods. Another extension could be considering multiple products into the proposed framework.
\section*{Appendix}
\appendix
\section{Linearization}\label{linearization}
\unskip
The proposed model in Section \ref{def} is non-linear due to the following situations:

\subsection{Multiplication of integer (binary) and continuous variables}
\begin{proposition} \label{prop:yint-relax}
Consider the term $q_{k}\cdot n_{k}$ in objective (\ref{obj1}). For bounded continuous and integer variables $q$ with $q_k \leq a_k , a_k > 0, \forall k \in K$ and $n$ with $n_k\leq b_k, b_k\geq 1, \forall k \in K$ , a non-negative bi-linear variables $\varsigma_k$ is defined as follows:
 $$\varsigma_k= q_{k}\cdot n_{k} \quad \forall k \in K $$
 with the following constraints added to the constraint set:
 \begin{alignat}{3}
& \varsigma_k \geq 0  \quad \forall k \in K &&\label{fdd}\\
& \varsigma_k \leq a_k q_k  \quad \forall k \in K &&\label{fdd1}\\
& \varsigma_k \leq b_k n_k  \quad \forall k \in K &&\label{fdd2}\\
& \varsigma_k \geq a_k q_k + b_k n_k - a_k b_k \quad \forall k \in K &&\label{fdd3}
\end{alignat}
For other variables with the same nature, we apply the same approach. 
\begin{proof}
For the proof, see \cite{gupte2013solving}.
\end{proof}
\end{proposition}
\subsection{Multiplication of multiple binary variables}
\begin{proposition} \label{prop:yint-relax_2}
Suppose that there are $n$ binary variables $x_i$ and the goal is to linearize the product of the variables as follows:
$$z=\prod_{i=1}^{n} x_i$$
where $z$ is a binary variable and denotes the value of the product. Then we have:
$$z\leq x_i \quad \forall i=1,...,n$$
$$z\geq \sum_{i=1}^{n}x_{i} - (n-1)$$
\begin{proof}
For the proof, see \cite{glover1973further}.
\end{proof}
Consider the term $t_k\cdot t^{'}_{k}$ in the objective (\ref{obj1}) where two binary variables are multiplied together. We define binary variables $\vartheta_{k}$ as follows:
$$\vartheta_{k}=t_k \cdot t^{'}_{k} \quad \forall k \in K$$
Then, the following constraints are added to constraints set:
\begin{alignat}{3}
& \vartheta_k \leq t_k  \quad \forall k \in K &&\label{fdd8}\\
& \vartheta_k \leq t^{'}_{k}  \quad \forall k \in K &&\label{fdd9}\\
& \vartheta_k \geq t_k+t^{'}_{k}-1  \quad \forall k \in K &&\label{fdd9}
\end{alignat}
\end{proposition}
Please note that we can always convert a bounded integer variable to the summation of some binary variables and use Proposition \ref{prop:yint-relax_2} for multiplication of integer and binary variables together.
\subsection{Square root of binary variables}
\begin{proposition} \label{prop:yint-relax_3}
Consider term $\sqrt{\sum_{i=1}^{n}a_{i}x_{i}}$ with $n$ binary variables $x_i$ and $a_i$s are parameters. We define parameter $b_j$ as follows:\\
$$b_j=\sum_{i,S_{i,j}=1} a_j \quad \forall j\in \{1,...,2^n\} $$
where $S$ is a binary matrix with $n$ rows and $2^n$ (number of subsets) columns. We define  binary variable $y_j$ where it is 1 if state $j$ is occurred with:
$$n y_{j}\leq \sum_{i=1,S_{i,j}=1}^{n} x_i +  \sum_{i=1,S_{i,j}=0}^{n} (1-x_i)\leq y_j +n -1 $$
then we can replace $\sqrt{\sum_{i=1}^{n}a_{i}x_{i}}$ with $\sum_{j=1}^{2^n}\sqrt{b_j}y_j$ 
\begin{proof}
For the proof, see \cite{boyd2004convex}.
\end{proof}
Consider term $\sqrt{\sum_{i \in I}\sigma^{2}y_{i,k}}$, following the above logic, we have :
$$b_{j,k}=\sum_{i\in I , S_{i,j,k}=1}\sigma_{i}^2 \quad \forall k \in K, j\in \{1,...,2^n\}$$
$$\sqrt{\sum_{i\in I}\sigma_{i}^{2}y_{i,k}}=\sum_{j=1}^{2^{|I|}}\sqrt{b_{j,k}}\tau_{j,k}$$
We define binary variable $\tau_{j,k}$ which is activated by the following constraints:
\begin{alignat}{3}
& \sum_{i,S_{i,j,k}=1}y_{i,k}+\sum_{i,S_{i,j,k}=0}(1-y_{i,k}) \geq |I|\cdot \tau_{j,k}  \quad \forall k \in K &&\label{fdd38}\\
&  \sum_{i,S_{i,j,k}=1}y_{i,k}+\sum_{i,S_{i,j,k}=0}(1-y_{i,k}) \leq \tau_{j,k} + |I|-1   \quad \forall k \in K &&\label{fdd28}
\end{alignat}
Following the three propositions, we linearized objective (1) and (2) as well as constraints (11)-(14),(22) and (23) to efficiently solve the problem.
\end{proposition}
% \subsection{}
% The appendix is an optional section that can contain details and data supplemental to the main text. For example, explanations of experimental details that would disrupt the flow of the main text, but nonetheless remain crucial to understanding and reproducing the research shown; figures of replicates for experiments of which representative data is shown in the main text can be added here if brief, or as Supplementary data. Mathematical proofs of results not central to the paper can be added as an appendix.

\section{Taguchi Experiments}\label{Taguhci}
% All appendix sections must be cited in the main text. In the appendixes, Figures, Tables, etc. should be labeled starting with `A', e.g., Figure A1, Figure A2, etc. 
The objective of parameter tuning is to calibrate the input parameters value of evolutionary algorithms. In this research, a three-level Taguchi design including 9 experiments for NSGA-II and NRGA and 27 experiments for SPEA-II and PESA-II is performed. The low, medium and high levels of the parameters along with their range are provided in Tables \ref{tab:parameters-1} and \ref{tab:parameters-2}.

\begin{table*}[!ht] 
	\caption{\hspace{2.2cm} Parameters for NSGA-II and NRGA}
	\centering
	\begin{tabular}{|c|c c c| c c c|} 
		\hline
		\multirow{1}{*}{Parameters}& \multicolumn{3}{c|}{NSGA-II}&\multicolumn{3}{c|}{NRGA}
		\\\cline{2-4} \cline{5-7}
		&  Low & Medium & High & Low & Medium &  High  \\
		\hline
		\multirow{1}{*}{Population Size}      & 50   & 100  & 150  & 50   & 100  & 150 \\
        \multirow{1}{*}{Crossover Percentage} & 0.70 & 0.80 & 0.90 & 0.70 & 0.80 & 0.90 \\
        \multirow{1}{*}{Mutation Percentage}  & 0.30 & 0.20 & 0.10 & 0.30 & 0.20 & 0.10 \\
        \multirow{1}{*}{Mutation Rate}        & 0.03 & 0.05 & 0.07 & 0.03 & 0.05 & 0.07 \\
		\hline
	\end{tabular}
	\label{tab:parameters-1}
\end{table*}

\begin{table*}[!htbp] 
	\caption{\hspace{2.2cm} Parameters for SPEA-II and PESA-II}
	\centering
	\begin{tabular}{|c|c c c| c c c|} 
		\hline
		\multirow{1}{*}{Parameters}& \multicolumn{3}{c|}{SPEA-II}& \multicolumn{3}{c|}{PESA-II}
		\\\cline{2-4} \cline{5-7}
		&  Low & Medium & High & Low & Medium &  High  \\
		\hline
		\multirow{1}{*}{Population Size}      & 50   & 100  & 150  & 50   & 100  & 150  \\
        \multirow{1}{*}{Archive Size}         & 100  & 200  & 300  & 100  & 200  & 300  \\
        \multirow{1}{*}{Crossover Percentage} & 0.70 & 0.80 & 0.90 & 0.70 & 0.80 & 0.90 \\
        \multirow{1}{*}{Mutation Percentage}  & 0.30 & 0.20 & 0.10 & 0.30 & 0.20 & 0.10 \\
        \multirow{1}{*}{Mutation Rate}        & 0.03 & 0.05 & 0.07 & 0.03 & 0.05 & 0.07 \\
        \multirow{1}{*}{Selection Pressure}   & -    & -    & -    & 2    & 3    & 4    \\
        \multirow{1}{*}{Deletion  Pressure}   & -    & -    & -    & 1    & 2    & 3    \\
		\hline
	\end{tabular}
	\label{tab:parameters-2}
\end{table*}

Each experiment is a different combination of factors specified where parameters change based on a predefined range. Each algorithm is tested in each experiment three times and the efficiency metrics are calculated. The average values for these runs are shown in Tables \ref{tab:adjust-param-1} and \ref{tab:adjust-param-2}.

%%%%%%%%%%%%%%%%%%%%%%%%%%%%%%%%%%%%
\begin{table*}
	\caption{Results from adjusting the parameters for  NSGA-II and NRGA}
	\centering
	\begin{tabular}{|c|c c c c|c c c c|} 
		\hline
		\multirow{1}{*}{NO.}& \multicolumn{4}{c|}{NSGA-II}& \multicolumn{4}{c|}{NRGA}
		\\\cline{2-5} \cline{6-9}
		&  QM & SM & MID & DM & QM & SM & MID & DM  \\
		\hline
		\multirow{1}{*}{1} & 0.96   & 168.03 &  0.80   &    193.54 &    0.80   &    147.97 &    0.79 & 186.67 \\
  		\multirow{1}{*}{2} & 0.69	& 107.57 &	0.79   &	177.33 &	0.77   &	146.32 &	0.79 & 184.48 \\
  		\multirow{1}{*}{3} &  0.84	& 163.82 &	0.79   &	188.95 &	0.73   &	133.59 &	0.78 & 182.52\\
        \multirow{1}{*}{4} &  0.98	& 176.86 &	0.79   &	195.40 &	0.87   &	150.96 &	0.79 & 191.21\\
        \multirow{1}{*}{5} &  0.84	& 132.69 &	0.80   &	185.02 &	0.73   &	155.65 &	0.81 & 188.12\\
        \multirow{1}{*}{6} &  0.98	& 179.64 &	0.79   &	195.48 &	0.84   &	150.47 &	0.80 & 187.48\\
        \multirow{1}{*}{7} &  0.93	& 159.19 &	0.80   &	191.67 &	0.87   &	168.84 &	0.80 & 194.29\\
        \multirow{1}{*}{8} &  0.87	& 159.41 &	0.80   &   	191.54 &	0.84   &	168.28 &	0.79 & 192.81\\
        \multirow{1}{*}{9} &  0.91	& 169.89 &	0.80   &	192.04 &	0.96   &	168.03 &	0.80 & 193.54\\
		\hline
	\end{tabular}
	\label{tab:adjust-param-1}
\end{table*}
 %%%%%%%%%%%%%%%%%%%%%%%%%%%%%%%%%%%%%%%%%%%%%%%
 \begin{table*} 
	\caption{Results from adjusting the parameters for  SPEA-II and PESA-II}
	\centering
	\begin{tabular}{|c|c c c c|c c c c|} 
		\hline
		\multirow{1}{*}{NO.}& \multicolumn{4}{c|}{SPEA-II}& \multicolumn{4}{c|}{PESA-II}
		\\\cline{2-5} \cline{6-9}
		&  QM & SM & MID & DM & QM & SM & MID & DM  \\
		\hline
		\multirow{1}{*}{1}  &   0.84  & 	138.29  &	0.80  &  185.55  &	0.64  &	101.05  &	0.78  &	 168.91\\
  		\multirow{1}{*}{2}  &   0.68  & 	112.08  &	0.79  &  176.22  &	0.62  &	77.42   &	0.79  &	 173.09\\
  		\multirow{1}{*}{3}  &   0.82  & 	156.65  &	0.80  &	 188.79  &	0.71  &	115.38  &	0.80  &	 180.04\\
        \multirow{1}{*}{4}  &   0.84  & 	161.96  &	0.80  &	 191.84  &	0.47  &	122.86  &	0.83  &	 189.39\\
        \multirow{1}{*}{5}  &   0.93  & 	166.22  &	0.79  &	 192.64  &	0.82  &	152.46  &	0.80  &	 189.91\\
        \multirow{1}{*}{6}  &   0.69  & 	149.58  &	0.80  &	 192.13  &	0.48  &	81.78   &	0.82  &	 161.79\\
        \multirow{1}{*}{7}  &   0.87  & 	146.77  &	0.79  &  186.84  &	0.41  &	68.92   &	0.78  &	 138.88\\
        \multirow{1}{*}{8}  &   0.81  & 	146.00  &	0.80  &	 188.93  &	0.44  &	102.36  &	0.79  &	 162.33\\
        \multirow{1}{*}{9}  &   0.84  & 	151.97  &	0.79  &	 189.42  &	0.43  &	69.03   &	0.83  &	 154.05\\
        \multirow{1}{*}{10} &   0.76  & 	128.16  &	0.80  &	 182.69  &	0.44  &	48.47   &	0.80  &	 145.41\\
        \multirow{1}{*}{11} &   0.81  & 	132.30  &	0.80  &	 180.82  &	0.51  &	53.41   &	0.79  &	 171.51\\
        \multirow{1}{*}{12} &   0.82  & 	155.13  &	0.80  &  188.37  &	0.52  &	91.80   &	0.78  &	 149.59\\
        \multirow{1}{*}{13} &   0.87  & 	156.95  &	0.80  &	 190.81  &	0.40  &	90.66   &	0.83  &	 166.53\\
        \multirow{1}{*}{14} &   0.87  & 	157.50  &	0.79  &	 189.20  &	0.69  &	118.82  &	0.79  &	 181.36\\
        \multirow{1}{*}{15} &   0.84  & 	158.91  &	0.80  &	 194.76  &	0.72  &	109.28  &	0.80  &	 174.61\\
        \multirow{1}{*}{16} &   0.89  & 	154.30  &	0.80  &	 191.12  &	0.64  &	119.73  &	0.82  &	 179.62\\
        \multirow{1}{*}{17} &   0.84  & 	145.79  &	0.80  &	 187.45  &	0.63  &	83.460  &	0.78  &	 159.3\\
        \multirow{1}{*}{18} &   0.91  & 	161.72  &	0.80  &	 193.39  &	0.60  &	149.84  &	0.82  &	 204.1\\
        \multirow{1}{*}{19} &   0.78  & 	142.58  &	0.80  &	 188.10  &	0.50  &	81.93   &	0.79  &	 149.84\\
        \multirow{1}{*}{20} &   0.86  & 	154.53  &	0.80  &	 189.93  &	0.44  &	57.72   &	0.81  &	 160.63\\
        \multirow{1}{*}{21} &   0.84  & 	144.64  &	0.79  &	 184.78  &	0.54  &	128.17  &	0.80  &	 182.44\\
        \multirow{1}{*}{22} &   0.80  & 	148.62  &	0.80  &	 188.84  &	0.64  &	118.59  &	0.80  &	 177.79\\
        \multirow{1}{*}{23} &   0.77  & 	145.81  &	0.80  &	 185.54  &	0.62  &	100.56  &	0.78  &	 172.19\\
        \multirow{1}{*}{24} &   0.87  & 	156.38  &	0.80  &	 191.46  &	0.42  &	88.44   &	0.82  &	 178.65\\
        \multirow{1}{*}{25} &   0.91  & 	156.46  &	0.79  &	 190.69  &	0.56  &	93.29   &	0.80  &	 174.68\\
        \multirow{1}{*}{26} &   0.91  & 	164.51  &	0.80  &	 193.37  &	0.48  &	67.00   &	0.81  &	 160.09\\
        \multirow{1}{*}{27} &   0.80  & 	147.19  &	0.80  &	 191.86  &	0.48  &	88.29   &	0.80  &	 160.05\\
		\hline
	\end{tabular}
	\label{tab:adjust-param-2}
\end{table*}
%%%%%%%%%%%%%%%%%%%%%%%%%%%%%%%%%%%%%%%%%%%%%%
\begin{figure}
\centering
\includegraphics[scale=0.6]{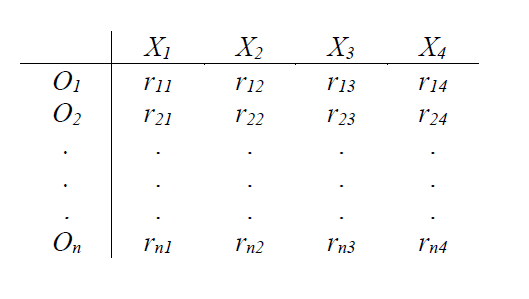}
\captionsetup{justification=centering}
\caption{\hspace{0.2cm} Experiments' result matrix}
\label{mat}
\end{figure}

We aim to transform each metric into a response. It should be noted that the higher values for QM and DM and lower values for MID and SM are preferable. Figure \ref{mat} is the matirx of results where  $X_j$ and $O_i$ denote $j-th$ metric and $i-th$ experiment, respectively, and $r_{ij}$ is the corresponding value. Since the metrics have different interpretation, we use a normalization method to convert all the values to comparable values. $X^{+}_j$ and $X^{-}_j$ represent normalized values for metrics QM/DM and MID/SM, respectively.
\begin{equation}
\stackrel{\text{+}}{X_j}\rightarrow R_{ij}=\frac{r_{ij} - \displaystyle\min_{i=1:n}(r_{ij})}{\displaystyle\max_{i=1:n}(r_{ij})-\displaystyle\min_{i=1:n}(r_{ij})}, \quad     \stackrel{\text{-}}{X_j}\rightarrow R_{ij}=\frac{\displaystyle\max_{i=1:n}(r_{ij})-r_{ij}}{\displaystyle\max_{i=1:n}(r_{ij})-\displaystyle\min_{i=1:n}(r_{ij})} 
\end{equation}

Using Goal Programming, we rank the metrics based on their importance and consider a weight for each of them and calculate the weighted sum (Response) as follows:

\begin{equation}
    Response_i = \sum_{j=1}^{4}R_{ij}w_j ;\quad (w_{QM}, w_{MID}, w_{SM}, w_{DM})=(10^2, 10, 1, 1)
\end{equation}

The weights are defined based on the importance of each metric. Normalized values and associated Responses for each algorithm are shown in Tables \ref{tab:responses-1} and \ref{tab:responses-2}.

%%%%%%%%%%%%%%%%%%%%%%%%%%%%%%%%%%%%%%%%%%%%%%
\begin{table*}[!htbp] 
	\caption{Results from adjusting the parameters for  NSGA-II and NRGA}
	\centering
	\begin{tabular}{|c|c c c c c|c c c c c|} 
		\hline
		\multirow{1}{*}{NO.}& \multicolumn{5}{c|}{NSGA-II}& \multicolumn{5}{c|}{NRGA}
		\\\cline{2-6} \cline{7-11}
		&  QM & SM & MID & DM & Response & QM & SM & MID & DM & Response \\
		\hline
		\multirow{1}{*}{1}  &   0.92	&	0.16	&	0.37	&	0.89	&	97.04	&	0.3	&	0.59	&	0.74	&	0.35	&	38.37\\
  		\multirow{1}{*}{2}  &   0	&	1	&	1	&	0	&	11	&	0.17	&	0.64	&	0.62	&	0.17	&	24.17\\
  		\multirow{1}{*}{3}  &   0.54	&	0.22	&	0.83	&	0.64	&	63.05	&	0	&	1	&	1	&	0	&	11\\
        \multirow{1}{*}{4}  &   1	&	0.04	&	0.69	&	1	&	107.9	&	0.6	&	0.51	&	0.45	&	0.74	&	65.76\\
        \multirow{1}{*}{5}  &   0.54	&	0.65	&	0.15	&	0.42	&	56.46	&	0	&	0.37	&	0	&	0.48	&	0.85\\
        \multirow{1}{*}{6}  &   1	&	0	&	0.69	&	1	&	107.94	&	0.5	&	0.52	&	0.41	&	0.42	&	55.05\\
        \multirow{1}{*}{7}  &   0.85	&	0.28	&	0.05	&	0.79	&	86.17	&	0.6	&	0	&	0.31	&	1	&	64.1\\
        \multirow{1}{*}{8}  &   0.62	&	0.28	&	0	&	0.78	&	62.6	&	0.5	&	0.02	&	0.73	&	0.87	&	58.18\\
        \multirow{1}{*}{9}  &   0.77	&	0.14	&	0.37	&	0.81	&	81.56	&	1	&	0.02	&	0.42	&	0.94	&	105.14\\
		\hline
	\end{tabular}
	\label{tab:responses-1}
\end{table*} 

%%%%%%%%%%%%%%%%%%%%%%%%%%%%%%%%%%%%%%%%%%%%%%%
\begin{table*}[!htbp] 
	\caption{Results from adjusting the parameters for  SPEA-II and PESA-II}
	\centering
	\begin{tabular}{|c|c c c c c|c c c c c|} 
		\hline
		\multirow{1}{*}{NO.}& \multicolumn{5}{c|}{SPEA-II}& \multicolumn{5}{c|}{PESA-II}
		\\\cline{2-6} \cline{7-11}
		&  QM & SM & MID & DM & Response & QM & SM & MID & DM & Response \\
		\hline
		\multirow{1}{*}{1}  &   0.64  & 0.52  &	0.16  &	0.5   &	66.28  &	0.57  &	0.49  &	0.9	  & 0.46  &	67.38\\
  		\multirow{1}{*}{2}  &   0	  & 1	  & 0.89  &	0	  & 9.86   &	0.52  &	0.72  &	0.71  &	0.52  &	60.39\\
  		\multirow{1}{*}{3}  &   0.55  &	0.18  &	0.41  &	0.68  &	59.72  &	0.73  &	0.36  &	0.58  &	0.63  &	80.13\\
        \multirow{1}{*}{4}  &   0.66  &	0.08  &	0.42  &	0.84  &	71.37  &	0.16  &	0.28  &	0	  & 0.77  &	17.03\\
        \multirow{1}{*}{5}  &   1	  & 0	  & 0.69  &	0.89  &	107.81 &	1	  & 0	  & 0.63  &	0.78  &	107.09\\
        \multirow{1}{*}{6}  &   0.04  &	0.31  &	0.43  &	0.86  &	9.29   &	0.19  &	0.68  &	0.31  &	0.35  &	22.72\\
        \multirow{1}{*}{7}  &   0.75  &	0.36  &	1	  & 0.57  &	86.09  &	0.02  &	0.8	  & 1	  & 0	  & 12.7\\
        \multirow{1}{*}{8}  &   0.54  &	0.37  &	0.1	  & 0.69  &	56.18  &	0.08  &	0.48  &	0.73  &	0.36  &	15.74\\
        \multirow{1}{*}{9}  &   0.66  &	0.26  &	0.84  &	0.71  &	75.65  &	0.05  &	0.8	  & 0.1	  & 0.23  &	6.96\\
        \multirow{1}{*}{10} &   0.31  &	0.7	  & 0.25  &	0.35  &	34.08  &	0.08  &	1	  & 0.53  &	0.1	  & 14.75\\
        \multirow{1}{*}{11} &   0.54  &	0.63  &	0.44  &	0.25  &	59.37  &	0.25  &	0.95  &	0.68  &	0.5	  & 33.71\\
        \multirow{1}{*}{12} &   0.55  &	0.2	  & 0.18  &	0.66  &	57.42  &	0.28  &	0.58  &	0.95  &	0.16  &	38.03\\
        \multirow{1}{*}{13} &   0.75  &	0.17  &	0	  & 0.79  &	76.12  &	0	  & 0.59  &	0.08  &	0.42  &	1.81\\
        \multirow{1}{*}{14} &   0.75  &	0.16  &	0.89  &	0.7	  & 84.93  &	0.68  &	0.32  &	0.73  &	0.65  &	76.36\\
        \multirow{1}{*}{15} &   0.66  &	0.14  &	0.08  &	1	  & 68.21  &	0.76  &	0.42  &	0.63  &	0.55  &	83.65\\
        \multirow{1}{*}{16} &   0.84  &	0.22  &	0.09  &	0.8	  & 86.02  &	0.57  &	0.31  &	0.3	  & 0.62  &	61.32\\
        \multirow{1}{*}{17} &   0.66  &	0.38  &	0.44  &	0.61  &	71.58  &	0.53  &	0.66  &	0.85  &	0.31  &	62.33\\
        \multirow{1}{*}{18} &   0.93  &	0.08  &	0.56  &	0.93  &	99.63  &	0.47  &	0.03  &	0.17  &	1	  & 49.49\\
        \multirow{1}{*}{19} &   0.39  &	0.44  &	0.33  &	0.64  &	43.88  &	0.22  &	0.68  &	0.81  &	0.17  &	31.04\\
        \multirow{1}{*}{20} &   0.73  &	0.22  &	0.57  &	0.74  &	79.88  &	0.1	  & 0.91  &	0.49  &	0.33  &	15.65\\
        \multirow{1}{*}{21} &   0.64  &	0.4	  & 0.91  &	0.46  &	73.67  &	0.33  &	0.23  &	0.64  &	0.67  &	40.4\\
        \multirow{1}{*}{22} &   0.48  &	0.33  &	0.22  &	0.68  &	51.64  &	0.56  &	0.33  &	0.6	  & 0.6	  & 62.4\\
        \multirow{1}{*}{23} &   0.36  &	0.38  &	0.6	  & 0.5	  & 42.58  &	0.52  &	0.5	  & 0.85  &	0.51  &	61.6\\
        \multirow{1}{*}{24} &   0.75  &	0.18  &	0.27  &	0.82  &	78.84  &	0.04  &	0.62  &	0.17  &	0.61  &	7.14\\
        \multirow{1}{*}{25} &   0.93  &	0.18  &	0.88  &	0.78  &	102.8  &	0.37  &	0.57  &	0.51  &	0.55  &	42.74\\
        \multirow{1}{*}{26} &   0.93  &	0.03  &	0.41  &	0.92  &	98.05  &	0.17  &	0.82  &	0.45  &	0.33  &	22.74\\
        \multirow{1}{*}{27} &   0.48  &	0.35  &	0.59  &	0.84  &	55.5   &	0.17  &	0.62  &	0.57  &	0.32  &	24.17\\
		\hline
	\end{tabular}
	\label{tab:responses-2}
\end{table*} 
%%%%%%%%%%%%%%%%%%%%%%%%%%%%%%%%%%%%%%%%%%%%%%%%%%%%

\begin{figure*}
\centering
	\includegraphics[scale=1.5]{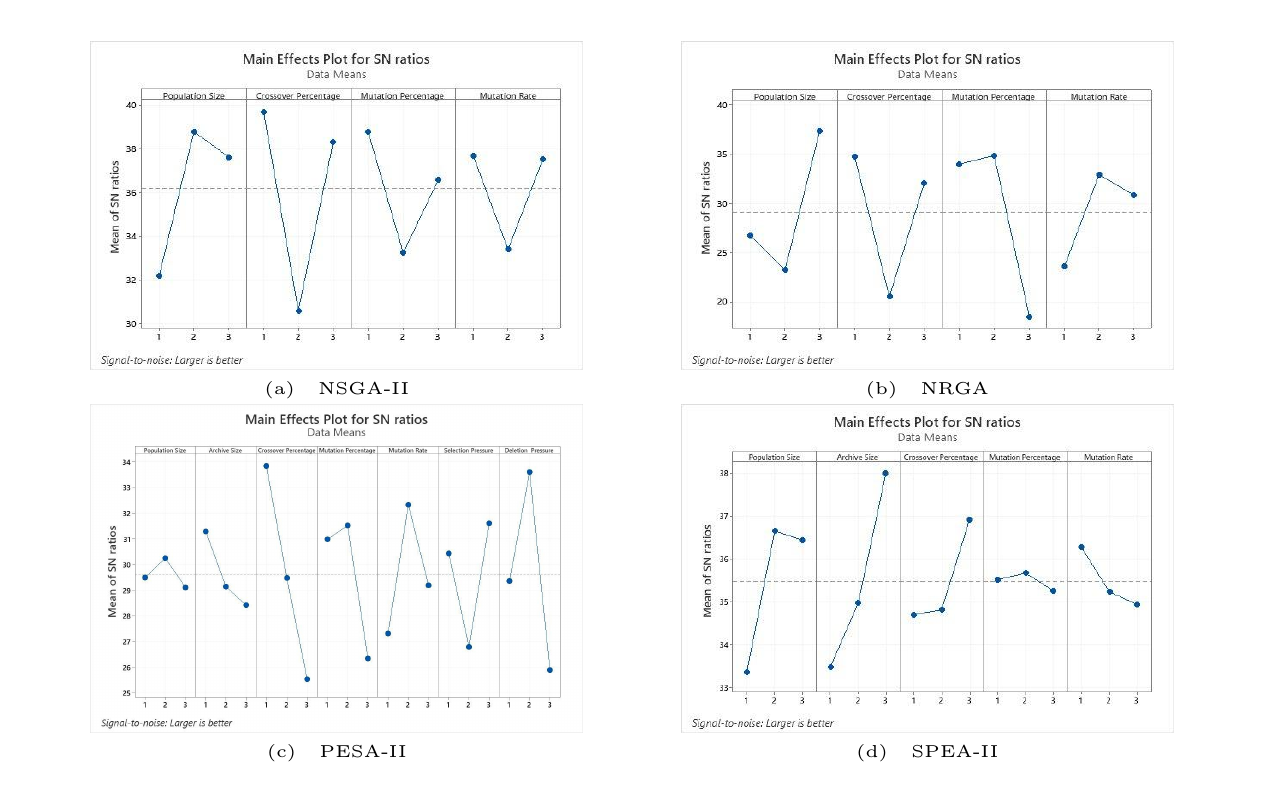}
	\captionsetup{justification=centering}
	\caption{S/N rates for each algorithm}
	\label{taguuuuchi}
\end{figure*}

Now, based on Response values, we calculate the S/N rates and determine the level of each input parameter as it is shown in Figure\ref{taguuuuchi}. The final values for adjusted parameters is provided in Table \ref{tuned}.

\newpage
\clearpage
\bibliographystyle{unsrt}
\bibliography{Refs}

%=====================================
% References, variant B: internal bibliography
%=====================================
% \begin{thebibliography}{999}
% % Reference 1
% \bibitem[Author1(year)]{ref-journal}
% Author1, T. The title of the cited article. {\em Journal Abbreviation} {\bf 2008}, {\em 10}, 142--149.
% % Reference 2
% \bibitem[Author2(year)]{ref-book}
% Author2, L. The title of the cited contribution. In {\em The Book Title}; Editor1, F., Editor2, A., Eds.; Publishing House: City, Country, 2007; pp. 32--58.
% \end{thebibliography}

% The following MDPI journals use author-date citation: Arts, Econometrics, Economies, Genealogy, Humanities, IJFS, JRFM, Laws, Religions, Risks, Social Sciences. For those journals, please follow the formatting guidelines on http://www.mdpi.com/authors/references
% To cite two works by the same author: \citeauthor{ref-journal-1a} (\citeyear{ref-journal-1a}, \citeyear{ref-journal-1b}). This produces: Whittaker (1967, 1975)
% To cite two works by the same author with specific pages: \citeauthor{ref-journal-3a} (\citeyear{ref-journal-3a}, p. 328; \citeyear{ref-journal-3b}, p.475). This produces: Wong (1999, p. 328; 2000, p. 475)

%%%%%%%%%%%%%%%%%%%%%%%%%%%%%%%%%%%%%%%%%%
%% Optional
% \sampleavailability{Samples of the compounds ...... are available from the authors.}

%% for journal Sci
%\reviewreports{\\
%Reviewer 1 comments and authors’ response\\
%Reviewer 2 comments and authors’ response\\
%Reviewer 3 comments and authors’ response
%}

%%%%%%%%%%%%%%%%%%%%%%%%%%%%%%%%%%%%%%%%%%
\end{document}